\documentclass[11pt]{amsart}
\usepackage{mathrsfs}
\usepackage{amssymb,latexsym}
\usepackage{amssymb,amsmath,ams}
\newdimen\AAdi%
\newbox\AAbo%
\def\AAk#1#2{\s_etbox\AAbo=\hbox{#2}\AAdi=\wd\AAbo\kern#1\AAdi{}}%
\def\AAr#1#2#3{\s_etbox\AAbo=\hbox{#2}\AAdi=\ht\AAbo\raise#1\AAdi\hbox{#3}}%
\font\tenmsb=msbm10 at 11pt \font\sevenmsb=msbm7 at 8pt
\font\fivemsb=msbm5 at 6pt
\newfam\msbfam
\textfont\msbfam=\tenmsb \scriptfont\msbfam=\sevenmsb
\scriptscriptfont\msbfam=\fivemsb

\newtheorem{theorem}{Theorem}[section]
\newtheorem{lemma}[theorem]{Lemma}
\newtheorem{definition}[theorem]{Definition}

\newtheorem{proposition}[theorem]{Proposition}

\newtheorem{corollary}[theorem]{Corollary}

\theoremstyle{remark}

\numberwithin{equation}{section}

\def\R{\mathbb R}
\def\S{\mathbb S}
\def\O{\mathcal O}
\def\a {\alpha}
\def\b{\beta}
\def\n{\nabla}
\def\g {\gamma}
\def\e {\eta}
\def\d {\delta}
\def\p {\partial}
\def\t {\theta}
\def\s{\mathscr{S}}
\def\k{\kappa}
\begin{document}
\title  [Convexity estimates of level sets]
{convexity estimates for level sets
of quasiconcave solutions to fully nonlinear elliptic equations}

\author{Pengfei Guan}
\address{Department of Mathematics\\
         McGill University\\
         Montreal, Quebec. H3A 2K6, Canada.}
\email{guan@math.mcgill.ca}
\author{Lu XU}
\address{Wuhan Institute of Physics and Mathematics\\
The Chinese Academy of Science\\
Wuhan,430071, HuBei Province, China}
\email{xulu@wipm.ac.cn}

\thanks{ Research of the first author was supported in part
by NSERC Discovery Grant. Research of the second author was
supported in part by NSFC No.10901159 and NSFC No.11071245.}

\begin{abstract} We establish a global geometric lower bound for the second fundamental form of the level surfaces of solutions to $F(D^2u, Du, u, x)=0$ in convex ring domains, in terms of boundary geometry and the structure of the elliptic operator $F$. We also prove a microscopic constant rank theorem,  
under a general structural condition introduced by Bianchini-Longinetti-Salani in \cite{BLS}.
\end{abstract}
 \maketitle

\section{Introduction}

Solutions of boundary value problems for elliptic equations often
inherit important geometric properties of the domains with the
influence of the structures of the corresponding equations. One of
these geometric features is the {\it quasiconcavity}. A function $u$
is called quasiconcave if its level sets $\{x| u(x)\ge c\}$ are
convex.
By the work of Gabriel \cite{Ga57}, the Green function of a
convex domain is quasiconcave. The same is also true for
$p$-harmonic functions in convex ring domains with homogeneous
boundary conditions following Lewis \cite{Le77}. Another example is
the quasiconcavity of solutions to the free boundary problem arising
in plasma physics in convex domains in the work of Caffarelli-Spruck
\cite{CS03}. The quasiconcavity of solutions to nonlinear equations
has been studied extensively in the literature, we refer \cite{Ga57, Le77,
CS82, Ka85, Ko90, CS03, Gr06, CS06, LS, BLS, Xu08, BGMX} and references therein.
The techniques of quasiconcave envelopes have been refined by Colesanti-Salani \cite{CS03},
and more recently by Bianchini-Longinetti-Salani \cite{BLS} to prove quasiconcavity of solutions to
general degenerate elliptic fully
nonlinear equations in the form
\begin{equation}\label{1.1}
F(D^2u, Du, u, x)=0,
\end{equation} in convex ring domain $\Omega=\Omega_0\setminus \Omega_1$ (i.e. $\Omega_0\subset \subset \Omega_1$ are convex)
with the Dirichlet boundary condition  \begin{eqnarray}\label{ringb}
u|_{\partial \Omega_0}=0, \quad \mbox{and} \quad  u|_{\partial
\Omega_1}=1.\end{eqnarray}

\medskip

The main focus of this paper is on the quantified properties of the quasiconcave solutions of equations of form (\ref{1.1}). More specifically, we establish a global a priori estimate on the geometric lower bound of the principal curvatures of the level surfaces of these quasiconcave solutions, in terms of boundary geometry and the structure of operator $F$. In addition to the geometric interest, this type of estimates may be used via homotopic deformation to obtain the existence of quasiconcave solutions of the corresponding equations. We achieve this macroscopic geometric estimate through a microscopic {\it Constant Rank Theorem} for the smallest principal
curvatures of the level surfaces for quasiconcave solutions. A simple microscopic convexity principle for level surfaces of solutions of equations in form (\ref{1.1}) is obtained in Theorem \ref{th1.1}, under a general structural condition introduced in \cite{BLS} to cover a larger class of nonlinear equations. A more refined version for the smallest principal curvatures of the level surfaces is proved in the last section of the paper. The main result, Theorem \ref{th1.2}, is a consequence of this type of microscopic principle.

\medskip

Let us introduce some notation. Denote $\s_n$ the space of real
symmetric $n \times n$ matrices and let $\Upsilon\subset \s_n$ be an
open set.
\begin{definition}\label{def-bt}
$\forall \; \t \in \S^{n-1}$, denote $\t^{\bot}$ the linear subspace
in $\R^n$ which is orthogonal to $\t$. Define
$\mathscr{S}^{-}_n(\t)$ to be the class of $n \times n$ symmetric
real matrices which are negative definite on $\t^{\bot}$. Denote
$\mathscr{S}^{0-}_n(\t)$ the subclass of $\mathscr{S}^{-}_n(\t)$ of
matrices that have $\t$ as eigenvector with corresponding null
eigenvalue. For any $b\in \R^n$ with $t=\langle b,\t\rangle>0$,
define
\begin{eqnarray}\label{B-a00}
\mathscr{B}^{-}_{\t}(\Upsilon)=\Big\{\mathbf   B \in
\mathscr{S}_{n+1} \;{\Big |}\;
\mathbf   B=\Big(\begin{array}{cc} \widetilde{\mathbf   B}&b^{T}\\
b&\chi
\end{array}\Big)\; \text{with}\quad \widetilde{\mathbf   B} \in
\mathscr{S}^{0-}_n({\t})\cap \Upsilon, \chi \in \R
\Big\}.\end{eqnarray}
\end{definition}

Denote $ \mathbf   J= \left(I_n \left|\mathbf   0\right.\right)$
the $n \times (n+1)$ matrix, where $I_n$ is the $n \times n$
identity matrix and $\mathbf   0$ is the null vector in $\mathbb
R^n$. Suppose $F=F(r,p,u,x)$ is a $C^{2}$ function in
$\Upsilon\times \R^n \times \R \times \Omega $,  $\forall
(\theta,u)\in \mathbb S^{n-1}\times \mathbb R$ fixed, set
\begin{equation}\label{Gamma} \Gamma_F=\Big\{(\mathbf   B,x)\in \mathscr{B}^{-}_{\t}(\Upsilon)\times\Omega: F(t^{-1}\mathbf   J\mathbf   B^{-1}\mathbf   J^T, t^{-1}\theta,
u,x)\ge 0\Big\}.\end{equation}

The following was proved in \cite{BLS}.
\begin{theorem}\label{th-bls} {\bf [Bianchini-Longinetti-Salani]}
Suppose $F$ is proper, continuous, degenerate elliptic operator
which satisfies a viscosity comparison principle. Assume that for
each $(\theta,u) $ fixed,  the super-level set $\Gamma_F$ defined in
(\ref{Gamma}) is convex. If $u \in C^{2}(\Omega)\bigcap C(\bar
\Omega)$ with $|\nabla u|>0$ is an admissible classical solution of
equation (\ref{1.1}) satisfying the Dirichlet boundary value
(\ref{ringb}) in convex ring domain $\Omega$, then level set $\{x\in
\Omega| u(x)\ge c\}\cup \Omega_1$ is convex for each constant $0\le
c\le 1$.
\end{theorem}

The class of operators $F$ satisfying
conditions in Theorem \ref{th-bls} includes Laplace operators,
$p$-Laplace operators, the Pucci operator, and the mean curvature type equations of the form
\begin{equation}\label{may-24} \sum_{i,j=1}^na_{ij}(\nabla u, u, x)u_{ij}=f(\nabla u,u,x). \end{equation}
A similar result was also proved by
Bianchini-Longinetti-Salani in \cite{BLS} under the assumption that, $\forall \theta, u$ fixed
\begin{equation}\label{cdA-0}  \Xi_{F} =\big\{(\mathbf   A,t,x)\in \Upsilon \times(0,+\infty)\times\Omega:
F(t^{-3}\mathbf   A, t^{-1}\theta, u,x)\ge 0\big\}\quad \mbox{is locally
convex.}\end{equation} With this structural condition on
$F$, a constant rank theorem was obtained in \cite{BGMX}. The
convexity structural condition on $F$ in Theorem \ref{th-bls}
is weaker than the convexity structural condition on $\Xi_{F}$. In particular, the mean curvature operator (\ref{may-24}) does not satisfy condition (\ref{cdA-0}).
Detailed discussion of these conditions as well as examples will be
given in section 2.

\medskip

To establish a strict convexity estimate on the second fundamental
forms of the level surfaces of solutions in Theorem \ref{th-bls}, we
need two assumptions:
\begin{equation}\label{1.2}
\mbox{\it ellipticity:} \quad \left(F^{\alpha\beta}\right):=\left(\frac{\partial
F}{\partial r_{\alpha\beta}}(\nabla^2 u(x), \nabla u(x), u(x),
x)\right) >0, \quad \forall x\in \Omega;
\end{equation}
\begin{equation}\label{1.4}
\mbox{\it structural condition:} \quad \quad \mbox{$\forall
(\theta,u) $ fixed,  the set $\Gamma_F$ is locally convex.}
\end{equation}
Throughout the paper, we assume
\begin{equation}\label{du} |Du(x)|\ge d_0>0, \forall \;x\in \Omega,\end{equation}
to ensure that the level-surface $\{x\in \Omega |
u(x)=c\}$ is smooth for each $c$.
\medskip

The first result of this paper is a microscopic constant rank theorem.
\begin{theorem}\label{th1.1}
Suppose $u \in C^{3,1}(\Omega)$ is a solution of (\ref{1.1}) and
$(D^2u(x), Du(x), u(x))\in \Upsilon\times \mathbb R^n\times
(-\gamma_0+\delta_0, \gamma_0+\delta_0)$ for some $\delta_0\in
\mathbb R$ at $x\in \Omega$. Suppose that $F$ satisfies conditions
(\ref{1.2}) -(\ref{1.4}) and the level set $\{x \in \Omega |u(x) \ge
c\}\cup \Omega_1$ of $u$ is connected and locally convex for all
$c\in (-\gamma_0+\delta_0, \gamma_0+\delta_0)$ for some
$\gamma_0>0$. Then the second fundamental form of level surface
$\Sigma^c=\{x \in \Omega |u(x) = c\}$ has the same constant rank for
all $c\in (-\gamma_0+\delta_0, \gamma_0+\delta_0)$.
\end{theorem}

\medskip

We now switch our attention to global geometric bounds of the second
fundamental forms of level surfaces of $u$. For a function $u$
defined in domain $\Omega$, denote 
\[\Sigma^{c}=\{x \in
\bar{\Omega} | u(x) = c\}\]
to be the level surface. For any  $x\in \Sigma^{c}$, denote $\k_s(x)$ the smallest
principal curvature of the level surface $\Sigma^{c}$ at $x$. For
each $c\in \mathbb R$, if $\Sigma^c\neq \emptyset$, set
\[\k^c=\inf_{x\in \Sigma^c} \k_s(x).\]
We will strengthen (\ref{1.2}) to
\begin{equation}\label{1.2-0}
\mbox{\it uniform ellipticity:} \quad \exists \lambda>0,
\big(\frac{\partial F}{\partial r_{\alpha\beta}}(D^2u(x), Du(x),
u(x), x)\big)\ge \lambda (\delta_{\a\b}) , \quad \forall  x\in
\bar\Omega.
\end{equation}
Set \begin{equation}\label{varpi} \varpi=\max_{\a,\b, \g, \e, i,j}\frac{\sup_{x\in \Omega}
\{|\frac{\partial^2 F(D^2u(x), Du(x), u(x), x)}{\partial
r_{\a\b}\partial r_{\g\e}}||\frac{\partial F(D^2u(x), Du(x), u(x),
x)}{\partial r_{ij}}||D u(x)|\}}{\lambda}.\end{equation}

\begin{theorem}\label{th1.2} Suppose $u$ is a classical solution of equation (\ref{1.1}) with the Dirichlet
boundary value (\ref{ringb}) in convex ring domain $\Omega$. Suppose
$F$ satisfies conditions
(\ref{1.4})-(\ref{1.2-0}) at $(D^2 u, Du, u, x) \in \Upsilon\times
\mathbb R^{n}\times [0,1] \times \Omega$. Then
\begin{equation}\label{k-est} \k^c \ge
\min\{\k^0e^{Ac}, \k^1e^{A(c-1)}, \frac{\lambda e^{A(c-1)}}{100\varpi} \},
\quad \forall \;c\in [0,1],\end{equation} for some universal
constant $A\ge 0$ depending only on $||F||_{C^{2}}, n, \lambda, d_0,
\|u\|_{C^{3}}$.
\end{theorem}

\medskip

It should be pointed out that the convexity estimates carried out in
this paper are very sensitive to the structure of the corresponding equation. For
equations of the form (\ref{may-24}) with the Dirichlet boundary
condition (\ref{ringb}), the behavior of $f$ is crucial. For
instance, in the case of Laplace equation
\begin{equation}\label{may-24-0} \Delta u=f(u),\end{equation}
Theorems \ref{th1.1}-\ref{th1.2} are true when $f\ge 0$. In general, Theorem \ref{th1.1} does not hold if $f(u)<0$ in equation (\ref{may-24-0}), even for $f\equiv -1$.

\bigskip

The rest of the paper is organized as follows. In Section 2, we
discuss the structural conditions and prove two key lemmas: Lemma \ref{lamma4.2-0} and Lemma
\ref{lem3.5-new}. An auxiliary curvature
test function is analyzed in Section 3. The proof of Theorem
\ref{th1.1} - Theorem \ref{th1.2} is given in the last
section, by establishing a strong maximum principle for the test function considered in Section 3.

\section{structural conditions}

We recall some notation and results in \cite{BLS}.
\begin{definition}
$\forall \; \t \in \S^{n-1}$, denote by
$\mathscr{A}^{-}_{\t}(\Upsilon)$ the following open set in
$\mathscr{S}_{n+1}$:
\begin{equation}\label{A-b0}
\mathscr{A}^{-}_{\t}(\Upsilon)=\left\{\mathbf   A \in
\mathscr{S}_{n+1}: \quad\right.
\mathbf   A=\Big(\begin{array}{cc} \widetilde{\mathbf   A}&\mu \t^{T}\\
\mu\t&0
\end{array}\Big) \quad \text{with} \quad \widetilde{\mathbf   A} \in
\left.\mathscr{S}^{-}_n(\t)\cap \Upsilon, \; \mu>0\right\}
\end{equation}\end{definition}

Properties of $\mathscr{A}^-_{\t}$, $\mathscr{B}^-_{\t}$ and their
relationship have been studied in \cite{BLS}. We list some of them
which will be used in this paper.

$det \mathbf   A \neq 0$ if $\mathbf   A \in \mathscr{A}^-_{\t}$,
and
\begin{equation}\label{B-a0}
\mathscr{B}^{-}_{\t}(\Upsilon)=\{\mathbf   A^{-1}:\; \mathbf   A \in
\mathscr{A}^-_{\t}(\Upsilon)\}.\end{equation}

If $\mathbf   B=\mathbf   A^{-1} \in \mathscr{B}^-_{\t}(\Upsilon)$,
then
\begin{equation}
\widetilde{\mathbf   A}=\mathbf   J\mathbf   B^{-1}\mathbf   J^T
\quad \text{and}\quad \mu={\frac1t}.
\end{equation}

Set
\begin{equation}\label{Q}
Q=t^2\mathbf   J\mathbf   B^{-1}\mathbf   J^T,\end{equation} where
$\mathbf   B\in \mathscr{B}^{-}_{\t}(\Upsilon)$ and $t$ defined in
Definition \ref{def-bt}. By symmetry of $\mathbf B$, $t=\langle b,\t
\rangle =\sum_{l=1}^nB_{n+1\;l}\t_l$.

In what follows, we will use summation over repeated indices
$\a,\b,\g,$ $\e, k,$ $l, m,$ $r, s \in \{1,...,n\}$ and $c, d, e, f
\in \{1,...,n+1\}$ unless otherwise indicated.

\begin{lemma}\label{lemG}[lemma~3.11 in \cite{BLS}]
$Q$ defined in (\ref{Q}) is concave in $\mathbf   B\in
\mathscr{B}^{-}_{\t}(\Upsilon)$. Furthermore
\begin{equation}\label{I-b0} I:=F^{\a\b}\frac{\p^2Q_{\a\b}}{\p B_{cd}\p B_{ef}}X_{cd}X_{ef} \le
0, \end{equation} for any nonnegative definite $n\times n$ matrix
$(F^{\a\b})$ and any $(n+1) \times (n+1)$ symmetric matrix
$(X_{cd})$.
 \end{lemma}
\noindent {\bf Proof.} The concavity of $Q$ has been proved in
\cite{BLS}. For any nonnegative definite $n\times n$ matrix
$(F^{\a\b})$, there exist $\zeta_1, \cdots, \zeta_n\in \mathbb R^n$,
such that
\[(F^{\a\b})=\zeta_1 \zeta_1^T+\cdots +\zeta_n \zeta_n^T.\]
Therefore, $I\le 0$ follows directly from the concavity of $Q$. \qed

\medskip

For function $F(r,p,u,x),$  write $F^{\alpha\beta}=\frac{\partial
F}{\partial r_{\alpha\beta}}, F^{p_l}=\frac{\partial F}{\partial
p_l}, \cdots$ as derivatives of $F$ with respect to corresponding
arguments. For the level set $\Gamma_F$ defined in (\ref{Gamma}),
denote the tangent space of $\Gamma_F$ as
\begin{equation*}\label{TF-0} \mathcal{T}\Gamma_F=\{V=((X_{cd}),(Z_k))\in \s_{n+1} \times\R^n:
\langle V, \nabla_{(\mathbf   B,x)} F(t^{-1}\mathbf   J\mathbf
B^{-1}\mathbf   J^T, t^{-1}\theta, u, x)\rangle=0\}. \end{equation*}

Write
\[F(t^{-1}\mathbf   J\mathbf   B^{-1}\mathbf   J^T,t^{-1}\theta,u,x)=F(t^{-3}Q,p, u, x).\]
Condition (\ref{1.4}) is equivalent to the fact
\[V\n_{(\mathbf   B,x)}^2FV^T \le 0, \quad \forall\; V=((X_{cd}),(Z_k))\in\; \mathcal{T}\Gamma_F.\]

A straight computation yields,
\begin{eqnarray}\label{fdF}
\n_{\mathbf   B}F&=&\Big(\frac{F^{\a\b}}{t^{4}}(t\frac{\partial
Q_{\a\b}}{\partial B_{cd}}-3Q_{\a\b}\d_{n+1
c}\d_{ld}\t_l)-\frac{F^{p_s}}{t^{2}}\delta_{n+1
c}\delta_{ld}\theta_s\theta_l\Big),\\
\label{fdF1} \nabla_x F&=&(F^{x_1},...,F^{x_n}),
\end{eqnarray}
\begin{eqnarray}\label{sdF}
V\n_{(\mathbf
B,x)}^2FV^T&=&\frac{F^{\a\b,\g\e}}{t^{8}}(t\frac{\partial
Q_{\a\b}}{\partial
B_{cd}}X_{cd}-3Q_{\a\b}X_{n+1l}\theta_l)(t\frac{\partial
Q_{\g\e}}{\partial
B_{ef}}X_{ef}-3Q_{\g\e}X_{n+1s}\theta_s)\nonumber\\
&&-2\frac{F^{\a\b,p_l}}{t^6}\t_l(t\frac{\partial Q_{\a\b}}{\partial
B_{cd}}X_{cd}-3Q_{\a\b}X_{n+1r}\theta_r)X_{n+1s}\t_s\nonumber\\
&&+\frac{F^{\a\b}}{t^3}\frac{\p^2Q_{\a\b}}{\p B_{cd}\p
B_{ef}}X_{cd}X_{ef}-6\frac{F^{\a\b}Q_{\a\b}}{t^5}X_{n+1l}X_{n+1s}\t_l\t_s\nonumber\\
&&-6\frac{F^{\a\b}}{t^5}(t\frac{\partial Q_{\a\b}}{\partial
B_{cd}}X_{cd}-3Q_{\a\b}X_{n+1r}\theta_r)X_{n+1s}\t_s\nonumber\\
&&+2\frac{F^{p_s}\theta_s}{t^{3}}X_{n+1l}X_{n+1r}\theta_l\t_r
+\frac{F^{p_s,p_m}\theta_s\t_m}{t^{4}}X_{n+1l}X_{n+1r}\theta_l\t_r+F^{x_k,x_l}Z_kZ_l\nonumber\\
&&+2\frac{F^{\a\b,x_k}}{t^4}(t\frac{\partial Q_{\a\b}}{\partial
B_{cd}}X_{cd}-3Q_{\a\b}X_{n+1r}\theta_r)Z_k-2\frac{F^{u_l,x_k}\t_l}{t^2}X_{n+1s}\t_sZ_k
.\end{eqnarray}

\medskip

This expression suggests us to set
\begin{eqnarray}\label{conv} \widetilde{X}_{\a\b}
=t^{-4}(\sum_{c,d}t\frac{\partial Q_{\a\b}}{\partial
B_{cd}}X_{cd}-3Q_{\a\b}\sum_{l}X_{n+1l}\theta_l), \quad
\widetilde{Y}=-t^{-2}\sum_{s}X_{n+1 s}\t_s.\end{eqnarray}
For
\[V=\big((X_{cd}),(Z_k)\big),\; H(V,V)=V\n_{(\mathbf B,x)}^2FV^T,\]
$H(V,V)$ can be written as
\begin{eqnarray}\label{3.61} H(V,V)&=&F^{\alpha\beta,
\g\e}\widetilde{X}_{\alpha\beta}\widetilde{X}_{\g\e}
+2F^{\alpha\beta,p_l}\theta_l\widetilde{X}_{\alpha\beta}\widetilde{Y}
+2F^{\a\b,x_k}\widetilde{X}_{\a\b}Z_k\nonumber\\
&&+F^{p_l, p_s}\t_l\t_s \widetilde{Y}^2+2F^{p_l,
x_k}\theta_l\widetilde{Y}Z_k+F^{x_k,x_l}Z_kZ_l+2tF^{p_l}\theta_l\widetilde{Y}^2\nonumber\\
&&+6tF^{\alpha\beta}\widetilde{X}_{\alpha\beta}\widetilde{Y}
-6t^{-1}F^{\alpha\beta}Q_{\alpha\beta}\widetilde{Y}^2 +\frac{I}{t^3}
\end{eqnarray}
where Einstein summation convention is used and $I$ is defined in
(\ref{I-b0}). At this point, we have proved
\begin{lemma}\label{lem3.5} Condition \eqref{1.4} is equivalent to $H(V,V)\le 0, \; \forall
\; V=(({X}_{cd}),(Z_k))\in \mathcal{T}\Gamma_F$, where
$F^{\alpha\beta, rs},F^{\a\b,p_l}, etc.$ in (\ref{3.61}) are
evaluated at $(t^{-3}Q, t^{-1}\theta, u, x)$.
\end{lemma}

We may now compare condition (\ref{1.4}) and (\ref{cdA-0}), these are the two structural
conditions introduced in \cite{BLS} ({\it Condition
(3.10)} and {\it Condition (1.2)} there).
As already discussed
by Bianchini-Longinetti-Salani in \cite{BLS}, a variation of these
two conditions can be compared (Theorem 3.12 in \cite{BLS}). In fact, the following is true.

\begin{corollary}\label{new0} The condition \eqref{cdA-0} that $\Xi_{F}$ for each $\theta, u$
is locally convex implies condition \eqref{1.4}.\end{corollary}
\noindent{\bf Proof.} Lemma~4.1 in \cite{BGMX} states that condition
\eqref{cdA-0} implies $H(V,V) \le t^{-3}I$ where $H(V,V)$ is defined
in (\ref{3.61}) and $I$ is defined in (\ref{I-b0}) respectively. The
corollary follows directly from lemma~\ref{lemG} and
lemma~\ref{lem3.5}.\qed

\medskip

The quantity $I$ defined in (\ref{I-b0}) is a crucial term. We wish to compute this term explicitly, so it can be used in the proof of main theorems in the last section. For our purpose, we set
$\t=(0,\cdots, 0, 1)$. In this case,   $\mathbf   A$
and $\mathbf   B$ can be written as (see \cite{BLS})
\begin{eqnarray}\label{A}
\mathbf   A&=& \left(
\begin{array}{ccccc}
 &&   &\times&0\\
&a_{ij}&&\vdots & \vdots\\
& & &\times&0\\
\times&\cdots&\times&\times&\mu\\
0&\cdots &0 &\mu&0
\end{array}\right),
\end{eqnarray}
\begin{eqnarray}\label{Bb}
\mathbf   B&=& \left(
\begin{array}{ccccc}
 &&&0&\times\\
&a^{ij}&&\vdots & \vdots\\
& & &0&\times\\
0&\cdots &0&0 &t\\
\times&\cdots&\times&t&\chi
\end{array}\right)
\end{eqnarray}
where the $(n-1) \times (n-1)$ matrix $(a_{ij})$ is negative
definite and can be assumed diagonal, $(a^{ij})$ is the inverse
matrix of $(a_{ij})$, $t=B_{n+1,n}={\frac1\mu}>0$. The values at the
positions denoted by $\times$ which are not important in the
calculations.

\medskip

Note that $B_{ln}=B_{nl}\equiv 0,\;\forall\; l\le n$. We may as well
set
\[X_{ln}=X_{nl}=0,\;\forall\;l\le n.\]
Denote
\[B^{\a\b}=(B^{-1})_{\a\b}=A_{\a\b}, \quad T:=\{1,..., n-1\}.\]
We compute
\begin{eqnarray}\label{sdQ}
&&F^{\a\b}\frac{\p^2Q_{\a\b}}{\p B_{cd}\p B_{ef}}X_{cd}X_{ef}\\
&=&2t^2F^{\a\b}B^{\a e}B^{fc}B^{d\b}X_{cd}X_{ef} -4tF^{\a\b}B^{\a
c}B^{d\b}X_{cd}X_{n+1n}+2F^{\a\b}B^{\a\b}X_{n+1n}^2,\nonumber
\end{eqnarray}
by breaking summation into the following
three parts.

\noindent {\it Case 1.} $\a,\b \in T$. We can see if $c=n$, then $d$
must be $n+1$, as
 $B^{n+1 \b}=0$,
\begin{eqnarray}\label{QT}
&& 2\sum_{i \in T}t^2B^{\a \a}B^{ii}B^{\b\b}X_{\a i}X_{\b i}
-4tB^{\a \a}B^{\b\b}X_{\a\b}X_{n+1n}+2F^{\a\a}B^{\a\a}X_{n+1
n}^2\nonumber\\
&=&\sum_{i\in T}\frac{2}{B^{ii}}(tB^{ii}B^{\a\a}X_{i\a}-B^{i
\a}X_{n+1n})(tB^{ii}B^{\b\b}X_{i\b}-B^{i \b}X_{n+1n}).
\end{eqnarray}

\noindent {\it Case 2.} $\a=n, \b\in T$ or $\b=n, \a\in T$. As $
B^{n n+1}=\dfrac{1}{t}$,
\begin{eqnarray}\label{QTN}
&& 2\sum_{i \in T}t^2B^{\a \a}B^{ii}B^{\b\b}X_{\a i}X_{\b i}
-4tB^{\a \a}B^{\b\b}X_{\a\b}X_{n+1n}+2F^{\a\a}B^{\a\a}X_{n+1
n}^2\nonumber\\
&=&2t^2B^{n e}B^{fc}(B^{\b\b}X_{c\b}+B^{n
\b}X_{cn})X_{ef}\nonumber\\
&&-4tB^{nc}(B^{\b\b}X_{c\b}+B^{n\b}X_{cn})X_{n+1n}+2B^{n\b}X_{n+1
n}^2\nonumber\\
&=&2t^2B^{n e}(\sum_{c\neq n}B^{fc}B^{\b\b}X_{c\b}+B^{f n+1}B^{n
\b}X_{ n+1n})X_{ef}\nonumber\\
&&-4t(\sum_{c\neq n}B^{nc}B^{\b\b}X_{c\b}+B^{n
n+1}B^{n\b}X_{n+1n})X_{n+1n}+2B^{n\b}X_{n+1
n}^2\nonumber\\
&=&2t^2B^{n e}\sum_{c\neq n}B^{fc}B^{\b\b}X_{c\b}X_{ef}
-4t\sum_{c\neq n}B^{nc}B^{\b\b}X_{c\b}X_{n+1n} \nonumber
\\
&=&2t^2\sum_{c\neq n}(\sum_{f\neq n}B^{n
e}B^{fc}B^{\b\b}X_{c\b}X_{ef}+B^{n n+1}B^{nc}B^{\b \b}X_{c\b}X_{
n+1n})\nonumber\\&&-4t\sum_{c\neq
n}B^{nc}B^{\b\b}X_{c\b}X_{n+1n} \nonumber\\
&=&2t^2\sum_{e\neq n}\sum_{i \in T}B^{n
e}B^{ii}B^{\b\b}X_{i\b}X_{ie}-2t\sum_{c\neq
n}B^{nc}B^{\b\b}X_{c\b}X_{n+1n} \nonumber\\
&=&2\sum_{i\in
T}\frac{1}{B^{ii}}(tB^{ii}B^{\b\b}X_{i\b}-B^{i\b}X_{n+1n})(t\sum_{c
\neq n}B^{ii}B^{nc}X_{ci}).
\end{eqnarray}

\noindent {\it Case 3.} $\a=\b=n$.
\begin{eqnarray}
&&\frac{\p^2Q_{nn}}{\p B_{cd}\p B_{ef}}X_{cd}X_{ef}\nonumber\\
&=&4t^2\sum_{c,d \neq n}(B^{nn}B^{n+1 c}B^{dn}X_{cd}X_{n n+1}+B^{n
n+1}B^{nc}B^{dn}X_{cd}X_{n+1n})\nonumber\\&&+2t^2\sum_{c,d,e,f \neq
n}B^{ne}B^{fc}B^{dn}X_{cd}X_{ef} +2t^2B^{nn}(B^{n n+1}B^{n+1
n}X_{n n+1}X_{n+1n}\nonumber\\
&&+B^{n+1 n+1}B^{n n}X_{n+1 n}X_{n n+1}+B^{n+1 n}B^{n+1
n}X_{n n+1}^2+B^{n n+1 }B^{n n+1 }X_{ n+1n}^2)\nonumber\\
&&-4t\sum_{c,d\neq n}B^{nc}B^{dn}X_{cd}X_{n+1 n}-8tB^{nn}B^{n+1
n}X_{n+1 n}^2+2B^{nn}X_{n+1 n}^2.\nonumber
\end{eqnarray}
It follows from the facts that $ B^{n n+1}=\dfrac{1}{t}$ and $B^{n+1
\b}=0, \quad \forall \b \neq n$,
\begin{eqnarray}\label{QN}
F^{nn}\frac{\p^2Q_{nn}}{\p B_{cd}\p B_{ef}}X_{cd}X_{ef}
&=&2t^2\sum_{c,d,e,f \neq n}F^{nn}B^{ne}B^{fc}B^{dn}X_{cd}X_{ef}\nonumber\\
&=&2t^2\sum_{i \in T}\sum_{d,e \neq
n}F^{nn}B^{ne}B^{ii}B^{dn}X_{id}X_{ei}\nonumber\\
&=&2\sum_{i \in T}\frac{F^{nn}}{B^{ii}}(t\sum_{e \neq
n}B^{ne}B^{ii}X_{ei})^2.
\end{eqnarray}

Set
\begin{equation}\label{Y-d0} Y_{i\a}:=tB^{ii}B^{\a\a}X_{i\a}-B^{i\a}X_{n+1n},  \; \forall \;\a
\in T; \quad  Y_{in}:=tB^{ii}\sum_{c \neq
n}B^{nc}X_{ci}.\end{equation} Combining \eqref{QT},\eqref{QTN} and
\eqref{QN}, for $\t=(0,\cdots,0,1)$, $I$ in (\ref{I-b0}) can be
written as
\begin{eqnarray}\label{sdQZ}
I &=&2\sum_{i \in T} \frac{F^{\a\b}}{B^{ii}}Y_{i\a}Y_{i\b},
\end{eqnarray} where $Y_{i\a},\; \forall \;\a \in \{1,...,n\}$ is defined in (\ref{Y-d0}).
\medskip

We wish to express $I$ in terms of $\widetilde X_{\a\b}$ and $\mathbf   A$.
Recall $Q_{\a\b}=t^2B^{\a\b}$,
\begin{equation}
\frac{\p Q_{\a\b}}{\p B_{cd}}=\frac{\p (t^2B^{\a\b})}{\p
B_{cd}}=-t^2B^{\a c}B^{d\b}+2tB^{\a\b}\d_{n+1 c}\d_{nd},
\end{equation}

\begin{equation}
\sum_{c,d}t\frac{\partial Q_{\a\b}}{\partial
B_{cd}}X_{cd}-3Q_{\a\b}\sum_{l}X_{n+1l}\theta_l=-t^3\sum_{c,d}B^{\a
c}B^{d\b}X_{cd}-t^2B^{\a\b}X_{n+1 n}.
\end{equation}
By \eqref{conv}
\begin{equation}\label{newX}  t^2\widetilde X_{\a\b}
=-tB^{\a c}B^{d\b}X_{cd}-B^{\a\b}X_{n+1,n}, \quad 1\le \a, \b \le n.\end{equation}
Extending the definition of $\widetilde X$ as $(n+1)\times (n+1)$ symmetric matrix by setting
\begin{equation}\label{newX0}  t^2\widetilde X_{ef}
=-tB^{e c}B^{df}X_{cd}-B^{ef}X_{n+1,n}, \quad 1\le e, f \le n+1.\end{equation}
Since $X_{n\a}=0, \forall \;\a\le n$, and $B^{n+1,c}=A_{n+1,c}=0,
\forall \;c\neq n$,
\begin{equation}\label{newX1} \widetilde X_{n+1,c}=0, \forall \;c\neq n; \quad \widetilde X_{n+1,n}=-\frac{2}{t^3}X_{n+1,n}=\frac2t\widetilde Y.\end{equation}
In this setting, $X_{cd}$ can be recovered using the formula below,
\begin{equation}\label{newX2}  X_{cd}=
-tA^{ce}A^{fd}\widetilde X_{ef}+\frac{t^2A^{cd}}2\widetilde X_{n+1,n}, \quad 1\le c,d\le n+1.\end{equation}

From the relationship $\mathbf   B={\mathbf   A}^{-1}$, and the fact that $(n-1)\times (n-1)$-matrix $(A_{\a\b})$ is diagonal, we have
$B_{n+1,i}=-\frac{tA_{n,i}}{A_{ii}}$. Hence from (\ref{Y-d0}),
\begin{equation}\label{Y-d2}
Y_{i\a}=t^3A_{i\a}\widetilde X_{n+1,n}-t^2\widetilde X_{i\a}.
\end{equation}

It follows from (\ref{sdQZ}) and (\ref{Y-d2}),
\begin{lemma}\label{lamma4.2-0} For $\t=(0,\cdots,0,1)$, if $(A_{\a\b})$ is diagonal,
then $I$ in (\ref{I-b0}) can be written as
\begin{eqnarray}\label{sdQZ-0}
I &=&2\sum_{i \in T} \frac{F^{\a\b}}{A_{ii}}Y_{i\a}Y_{i\b},
\end{eqnarray} where $Y_{i\a}$ is defined in (\ref{Y-d2}).
\end{lemma}

\medskip

From (\ref{fdF}), (\ref{fdF1}) and (\ref{conv}),
\begin{eqnarray*}
<V,\nabla_{(B,x)}F>={F^{\a\b}}\widetilde{X}_{\a\b}+F^{u_n}\widetilde{Y}+F^{x_k}Z_k
\end{eqnarray*}
$V=((X_{cd}),(Z_k))\in \mathcal{T}\Gamma_F$ if and only if
\begin{equation}\label{Vb1}
{F^{\a\b}}\widetilde{X}_{\a\b}+F^{u_n}\widetilde{Y}+F^{x_k}Z_k=0,\end{equation}
where $\widetilde{X}_{\a\b}$, $\widetilde{Y}$ as in (\ref{newX}) and
(\ref{newX1}) respectively, and $F^{\a\b}, F^{u_n}$ and $F^{x_k}$
are evaluated at $(t^{-1}\widetilde {\mathbf   A}, t^{-1}\theta,
u,x)$ with $\theta=(0,\cdots, 0, 1)$.

Set $\widetilde V =((\widetilde X_{cd}),(Z_k))
$ where $(\widetilde X_{cd})$ defined by \eqref{newX0} with
\eqref{newX1}, rewrite \eqref{3.61} as
\begin{eqnarray}\label{HV-0} H(\widetilde V, \widetilde V)&=&\frac{I}{t^3}+S,\end{eqnarray}
where $I$ is defined as in (\ref{sdQZ-0}),
\begin{eqnarray}\label{defS}
S&=& F^{\alpha\beta,
\g\e}\widetilde{X}_{\alpha\beta}\widetilde{X}_{\g\e}+
2F^{\alpha\beta,p_n}\widetilde{X}_{\alpha\beta}\widetilde{Y}
+2F^{\a\b,x_k}\widetilde{X}_{\a\b}Z_k+F^{p_n, p_n}\widetilde{Y}^2\nonumber\\
&&+2F^{p_n,
x_k}\widetilde{Y}Z_k+F^{x_k,x_l}Z_kZ_l
+2tF^{p_n}\widetilde{Y}^2+6tF^{\alpha\beta}\widetilde{X}_{\alpha\beta}\widetilde{Y}
-6tF^{\alpha\beta}A_{\alpha\beta}\widetilde{Y}^2
,
\end{eqnarray}
and $F^{\a\b}, F^{u_n}$, $F^{x_k}$ etc. are
evaluated at $(t^{-1}\widetilde {\mathbf   A}, t^{-1}\theta, u,x)$,
$\widetilde{\mathbf   A} \in \mathscr{S}^{-}_n(\t)\cap \Upsilon$.

\medskip

By (\ref{newX2}) and Lemma \ref{lem3.5},
\begin{lemma}\label{lem3.5-new} Condition \eqref{1.4} is equivalent to
\[H(\widetilde V,\widetilde V)\le 0, \;\forall
\; \tilde V=((\widetilde X_{cd}),(Z_k))  \mbox{ satisfying
(\ref{Vb1})}.\] By approximation, if Condition \eqref{1.4} is
satisfied, then \[H(\widetilde V,\widetilde V)\le 0,\] for every
$\tilde V=((\widetilde X_{cd}),(Z_k))$ satisfying (\ref{Vb1}) at
each $\widetilde{\mathbf   A} \in
\overline{\mathscr{S}^{-}_n(\t)}\cap \Upsilon$ diagonal, where
$\overline{\mathscr{S}^{-}_n(\t)}$ is the closure of
$\mathscr{S}^{-}_n(\t)$, with $Y_{k\a}=0$ when $a_{kk}=0$ for some
$k\le n-1$.
\end{lemma}

\bigskip

With the explicit expression of $H$ in Lemma \ref{lem3.5-new}, we may verify condition \eqref{1.4} for mean curvature operator and general quasilinear operator $F$ satisfying structural conditions in \cite{Xu08}.
Condition (\ref{cdA-0}) is not satisfied by mean curvature operator as indicated in \cite{BLS}.
It was verified there that for $n=2$, the Mean Curvature operator:
\begin{equation}\label{mc}
F(D^2u, Du)=Div(\dfrac{Du}{\sqrt{1+|Du|^2}})=\dfrac{\Delta
u}{\sqrt{1+|Du|^2}}-\sum_{\a,\b=1}^n\dfrac{u_{\a}u_{\b}u_{\a\b}}{\sqrt{1+|Du|^2}^3}=f(u)\ge 0,
\end{equation}
satisfies condition \eqref{1.4}, but not \eqref{cdA-0}. Here we verify
this fact for general $n$. Since condition \eqref{1.4} and
\eqref{cdA-0} are invariant under orthogonal transformation, we may
as well set $Du=(0,...,0,u_n)$, $(u_{ij})$ is diagonal for each $i,j
\in T=\{1,...,n-1\}$. We also note $(u_{ij})$ is negative definite.
According to \eqref{Gamma}, $t^{-1}\t=Du$, where $\t=(0,...,0,1)$ and
$t^{-1}=u_n$, $\mathbf   JB^{-1}\mathbf   J^T=u_n^{-1}(D^2u)$.
Since the mean Curvature operator $F$ in \eqref{mc} is homogenous of one
degree, $S$ in (\ref{defS}) can be calculated as
\begin{eqnarray}\label{Sbbb}
S&=&2F^{\alpha\beta,u_n}\widetilde{X}_{\alpha\beta}\widetilde{Y}
+F^{u_n,
u_n}\widetilde{Y}^2+2\frac{F^{u_n}}{u_n}\widetilde{Y}^2+6\frac{F^{\alpha\beta}}{u_n}\widetilde{X}_{\alpha\beta}\widetilde{Y}
-6u_n^{-2}F\widetilde{Y}^2
\end{eqnarray}

From \eqref{Vb1}, $\widetilde V=((\widetilde X_{cd}))$ satisfies
\begin{equation}\label{nF1}0= \langle \widetilde{V}, \n_B
F\rangle=F^{\a\b}\widetilde{X}_{\a\b}+F^{u_n}\widetilde{Y}.
\end{equation}
A straightforward calculation yields that
\begin{eqnarray}\label{Smc}
S&=& \sum_{i \in
T}\left(\dfrac{4u_{n}}{W^3}\widetilde{X}_{ii}\widetilde{Y}-\frac{6}{W^3}u_{ii}\widetilde{Y}^2\right)
-\dfrac{3+6u_n^{-2}}{W^4}F\widetilde{Y}^2
\end{eqnarray}
where $W=\sqrt{1+|Du|^2}$. It is easy to check that $S\le 0$ is violated for some
$\widetilde{X}_{ii},\widetilde{Y}$ satisfying \eqref{nF1}. On the other hand, Lemma 4.1 in \cite{BGMX} implies that $S\le 0$ if condition \ref{cdA-0} is satisfied. Therefore, $F$ does not satisfy condition \eqref{cdA-0}.

However, from \eqref{sdQZ-0}
\begin{eqnarray}\label{Imc}I
&=&2u_n\sum_{i \in T} \frac{F^{\a\b}}{u_{ii}}Y_{i\a}Y_{i\b}.
\end{eqnarray}
For the mean curvature equation, it can be computed that
\begin{equation}
\label{Yia}Y_{i\a}= -t^2\widetilde{X}_{i\a} +
2t^2u_{i\a}\widetilde{Y}u_n^{-1}, \quad \forall \;i, \a \in
G.\end{equation} By \eqref{Smc}, \eqref{Imc} and \eqref{Yia}, and the facts that $u_{ii}<0$, $F \ge 0$,
\begin{eqnarray*}\label{S+I}
H&=&S+Iu_n^3 \nonumber\\
 &\le& S+2u_n^4\sum_{i \in T}
\frac{F^{ii}}{u_{ii}}Y_{ii}^2\nonumber\\
&=& 2\sum_{i \in T}\left(\dfrac{1}{Wu_{ii}}\widetilde{X}_{ii}^2-2\dfrac{u_{n}
+2u_n^{-1}}{W^3}\widetilde{X}_{ii}\widetilde{Y}+\frac{4u_n^{-2}+1}{W^3}u_{ii}\widetilde{Y}^2\right)
-\dfrac{3+6u_n^{-2}}{W^4}F\widetilde{Y}^2\nonumber\\
&=&\sum_{i \in
T}\dfrac{2}{Wu_{ii}}\left[\left(\widetilde{X}_{ii}^2-\dfrac{u_{n}
+2u_n^{-1}}{W^2}u_{ii}\widetilde{Y}\right)^2+\frac{1}{W^4}u^2_{ii}\widetilde{Y}^2\right]
-\dfrac{3+6u_n^{-2}}{W^4}F\widetilde{Y}^2\nonumber\\
&\le& 0,
\end{eqnarray*}
That is the mean curvature operator $F$ satisfies
condition \eqref{1.4} by Lemma~\ref{lem3.5-new}. This example indicates that the term $I$ is the key.
The verification of condition \eqref{1.4} for the quasilinear operators considered in \cite{Ko90, Xu08}
can be done in a similar way, we leave for the interested reader to check them.

\section{the test function}

The proof of our main results relies on the establishment of a maximum principle
for certain appropriate curvature test function. This section is devoted to discuss some regularity and concavity properties of the proposed test function.

\medskip

{\it We will assume $u \in C^{3,1}(\Omega)$, $|\nabla u|>0$ and $\{x\in\Omega|u(x)\ge
c\}\cup\Omega_1$ is locally convex in the rest of this paper}.

We recall some of formulas related to the Weingarten curvature tensor of level surfaces.
Suppose $u$ is a function defined in an open set in $\mathbb R^n$,
assume that $u_n(x) \neq 0$. The upward inner normal direction of
the level sets of $u$ is
\begin{eqnarray}\label{2.12}
\vec{n}= \frac{|u_n|}{|Du|u_n}(u_1,u_2,...,u_{n-1},u_n).
\end{eqnarray}
It's calculated in \cite{BGMX} that the second fundamental form $II$ of
the level surface of function $u$ with respect to the upward
normal direction (\ref{2.12}) is
\begin{equation}\label{2.17}
h_{ij} =  - \frac{|u_n|(u_n^2 u_{ij} + u_{nn}u_{i}u_j - u_nu_ju_{in}
- u_nu_iu_{jn})}{|Du|u_n^3}, \quad i,j \le n-1.
\end{equation}

Note that as $\{x\in\Omega|u(x)\ge c\}\cup\Omega_1$ is locally
convex, the second fundamental form of $\Sigma^{c}$ is
nonnegative definite with respect to the upward normal direction
(\ref{2.12}). For $Du$ is as the same direction as $\vec{n}$, thus
we have $u_n>0$ locally. (\ref{2.17}) implies that the
matrix $(u_{ij}(x))$ is nonpositive definite.

Denote $a(x)=(a_{ij}(x))$  the symmetric Weingarten tensor of
$\Sigma^{u(x)}=\{y \in \Omega| u(y) =u(x)\}$. Our assumption implies
that $a$ is
nonnegative definite. Since $u_n(x)\neq 0$, following \cite{CNS85},
the Weingarten tensor can be computed as (see \cite{BGMX}),
\begin{equation}\label{2.18}
a_{ij} =h_{ij} -\frac{\sum_{l=1}^{n-1}u_iu_lh_{jl}}{W(1+W)u_n^2}
-\frac{\sum_{l=1}^{n-1}u_ju_lh_{il}}{W(1+W)u_n^2} +
\frac{\sum_{l,k=1}^{n-1}u_iu_ju_ku_l h_{kl}}{W^2(1+W)^2u_n^4}, \quad
i,j\le n-1,
\end{equation} where $W =
(1+|\nabla_{x'}v|^2)^{\frac12}$ and $x'=(x_1,\cdots,x_{n-1})$.

\medskip

Set
\begin{equation}\label{a}
\widetilde{a}=a-\eta_0g(u)I, \quad \eta_0 \ge 0,\quad g(u)=e^{Au}
\end{equation}
where $\eta_0\ge 0$ and $A \ge 0$ are constants to be determined later such that $\widetilde{a}\ge 0$.

Suppose the minimal
rank $l$ of $\widetilde{a}$ is attained at some interior point $x_0$. Let $\O$ be a small open neighborhood of $x_0$ such that for each $x \in \O$,
there are $l$ "good" eigenvalues of $(\widetilde{a}_{ij})$ which are
bounded below by a positive constant, and the other $n-1-l$ "bad"
eigenvalues of $(\widetilde{a}_{ij})$ are very small. Denote $G$ be
the index set of these "good" eigenvalues and $B$ be the index set
of "bad" eigenvalues. For each $x \in \O$ fixed, we may express
$(a_{ij})$ in a form of (\ref{2.18}), by choosing $e_1,\cdots,
e_{n-1},e_n$ such that
\begin{eqnarray}\label{coord}
\quad \quad \quad \mbox{\it $|Du|(x)=u_n(x)>0$,
$(u_{ij}(x)),i,j=1,..,n-1$ is diagonal.}\end{eqnarray} From
\eqref{2.18} and \eqref{a}, the matrix
$(\widetilde{a}_{ij}),i,j=1,..,n-1$ is also diagonal at $x$, and
without loss of generality we may assume $\widetilde{a}_{11} \le
\widetilde{a}_{22} \le ... \le \widetilde{a}_{n-1, n-1}$. There is a
positive constant $C>0$ depending only on $\|u\|_{C^{4}}$ and $\O$,
such that $\widetilde{a}_{n-1,n-1} \ge \widetilde{a}_{n-2,n-2}
\ge...\ge \widetilde{a}_{n-l,n-l} > C$ for all $x \in \O$. For
convenience we denote $G=\{n-l,n-l+1,...,n-1\}$ and
$B=\{1,2,...,n-l-1\}$ be the "good" and "bad" sets of indices
respectively. If there is no confusion, we also denote
\begin{eqnarray}\label{G-B}
\quad \quad \mbox{\it
$B=\{\widetilde{a}_{11},...,\widetilde{a}_{n-l-1,n-l-1}\}$ and
$G=\{\widetilde{a}_{n-l,n-l},...,\widetilde{a}_{n-1,
n-1}\}$.}\end{eqnarray} Note that for any $\delta>0$, we may choose
$\O$ small enough such that $\widetilde{a}_{jj} <\delta$ for all $j
\in B$ and $x \in \O$.

The following two functions are of fundamental importance in our treatment.

\begin{eqnarray}\label{3.1}
 p(\widetilde{a})=\sigma_{l+1}(\widetilde{a}_{ij}), \quad
q(\widetilde{a}) &=& \left\{
 \begin{array}{llr}
\frac{\sigma_{l+2}(\widetilde{a}_{ij})}{\sigma_{l+1}(\widetilde{a}_{ij})},  &\text{if} \; \sigma_{l+1}(\widetilde{a}_{ij})>0&\\
0,& otherwise. &
\end{array}
 \right.
\end{eqnarray}
 We consider function
\begin{equation}
 \label{3.4}
 \phi(\widetilde{a})= p(\widetilde{a})+q(\widetilde{a})
\end{equation}
where $p$ and $q$ as in \eqref{3.1}. The function $\phi$ was first introduced in \cite{BG08} for the Hessian of solution $u$, and for Weingarten tensor $a$ in \cite{BGMX}. Here we adopt it as a function in $\widetilde{a}$.

We will use notion $h=O(f)$ if $|h(x)| \le Cf(x)$ for $x \in \O$ with
positive constant $C$ under control.
Again, as in \cite{BG08}, to get around $p=0$, for $\varepsilon>0$ sufficiently small,
we instead consider
\begin{equation}\label{3.5}
\phi_\varepsilon (\widetilde{a})= \phi(\widetilde{a}_\varepsilon),
\end{equation}
where $\widetilde{a}_\varepsilon=\widetilde{a}+\varepsilon I.$ We
will also denote $G_\varepsilon=\{\widetilde{a}_{ii}+\varepsilon,
i\in G\},$ $B_\varepsilon=\{\widetilde{a}_{ii}+\varepsilon, i\in
B\}.$

We will write $p$ for $p_\varepsilon,$ $\phi$ for $\phi
_\varepsilon$, $q$ for $q_\varepsilon$, $\widetilde{a}$ for
$\widetilde{a}_\varepsilon$, $G$ for $G_\varepsilon$, $B$ for
$B_\varepsilon$ with the understanding that all the estimates will
be independent of $\varepsilon.$ In this setting, if we pick $\O$
small enough,
 there is $C>0$ independent of $\varepsilon$ such that
 \begin{equation}\label{3.6}
 \phi(\widetilde{a}(z))\geq C\varepsilon, \quad \sigma_1(B(z))\geq C\varepsilon,
 ~\quad ~\rm{ for ~all~ }z\in \O.
 \end{equation}

In what follows, $i, j, \cdots$ will be denoted as indices run from $1$ to
$n-1$ and the Greek indices
$\alpha, \beta, \cdots$ will be denoted as indices run from $1$ to $n$.
Denote
\begin{eqnarray*}
p_{\alpha}=\frac{\partial p}{\partial x_{\alpha}},\quad
p_{\alpha\beta}=\frac{\partial^2
p}{\partial x_{\alpha}\partial x_{\beta}},\quad
F^{\alpha \beta}=\frac{\partial F}{\partial u_{\alpha\beta}}, \quad 1\le
\alpha,\beta \le n.
\end{eqnarray*}
We also denote $g=e^{Au}$,
\begin{equation} \mathcal H_{\phi} = \sum_{i,j\in B}|\nabla \widetilde{a}_{ij}|+\phi,
\end{equation}
and $\forall j\in B$,
\begin{eqnarray}\label{I-j}
I_j=\sum_{i\in G}[-2u_n^3\sum_{\alpha,\beta\notin B}
\dfrac{F^{\alpha\beta}{{a}_{ij,\alpha}{a}_{ij,\beta}}}{{a}_{ii}}
+4u_n^{2}u_{nj}\sum_{\a \notin B} F^{\alpha
i}{a}_{ij,\a}+2u_{nj}^2F^{ii}{u}_{ii}],
\end{eqnarray}
and
\begin{eqnarray}\label{Jp}
J_{1j}&=& -12\sum_{\a=1}^nF^{j\a}u_{n\a}u_{jn}u_{jj}
+4\sum_{\a=1}^nF^{j\a}u_{jn\a}u_{n}u_{jj}-2u_{nn}F^{jj}u_{jj}^2\nonumber\\
&&-\eta_0\sum_{\alpha,\beta=1}^n F^{\alpha\beta}g_{\a\b}u_n^3+4\eta_0\sum_{\a=1}^nF^{j\a}g_{\a}u_{jn}u_n^2-2\eta_0\sum_{\alpha,\beta=1}^n F^{\alpha\beta}u_{n\a}g_{\b}u_n\nonumber\\
&&-\eta_0g\sum_{\alpha,\beta=1}^nF^{\alpha\beta}(2u_{j\a}u_{j\b}u_n+u_{n\a\b}u_n^2
+\sum_{i=1}^{n-1}u_{i\a}u_{i\b}u_n)
-2\eta_0g(\sum_{i\in
B}F^{ii})u_{nj}^2u_n\nonumber\\
&&-2\eta_0g\sum_{\alpha,\beta=1}^n\sum_{i\in G}
\dfrac{F^{\alpha\beta}{a}_{ij,\alpha}{a}_{ij,\beta}}{{a}_{ii}\widetilde{a}_{ii}}u_n^3,
\end{eqnarray}
and
\begin{eqnarray}\label{J2}
J_{2j}&=&2F^{\a\b,u_j}u_{j\a\b}u_{jj}+2F^{u_nu_j}u_{jn}u_{jj}+2F^{u_j,x_j}u_{jj}+F^uu_{jj}\\
&&+F^{u_l}u_{ln}u_n^{-1}u_{jj} +2F^{u_j}u_{jn}u_n^{-1}u_{jj}
+F^{u_ju_j}u_{jj}^2-\eta_0F^{u_l}g_lu_n.\nonumber
\end{eqnarray}

\begin{lemma}\label{lem3.3} Suppose $u\in C^{3,1}$ is a solution of equation (\ref{1.1}) with $|\nabla u|>0$, then $\phi \in C^{1,1}(\O)$.
For any fixed $x \in \O$, with
the coordinate chosen as in (\ref{coord}) and (\ref{G-B}),
\begin{equation}\label{3.8}
\phi_\alpha=\big[\sigma_l({G})+\frac{{\sigma}^2_1(B|j)
-{\sigma}_2(B|j)}{{\sigma}^2_1(B)}\big]\widetilde{a}_{jj,\alpha}
+O(\mathcal H_{\phi}),
\end{equation}
and
\begin{eqnarray}\label{3.37}
&&F^{\alpha\beta}\phi_{\alpha\beta}\nonumber\\&=& \sum_{j\in
B}u_n^{-3}\left[\sigma_l(G)+\frac{{\sigma}^2_1(B|j)-{\sigma}_2(B|j)}{{\sigma}^2_1(B)}\right]\left\{
 \Big[\sum_{\alpha,\beta,\g,\e=1}^n
F^{\alpha\beta, \g\e}u_{\alpha\beta j}u_{\g\e j}\right.\nonumber\\
&&+2\sum_{\alpha,\beta=1}^nF^{\alpha\beta,
u_n}u_{j\alpha\beta}u_{jn}+2\sum_{\alpha,\beta=1}^nF^{\alpha\beta,
x_j}u_{j\alpha\beta} +F^{u_n,u_n}u_{jn}^2+2F^{u_n,x_j}u_{jn}+F^{x_j,
x_j}\nonumber\\
&&\left.+2\frac{F^{u_n}}{u_n}u_{jn}^2\Big]u_n^2+6\sum_{\alpha,\beta=1}^nF^{\alpha\beta}u_{j\alpha\beta}u_{jn}u_n
-6\sum_{\alpha,\beta=1}^nF^{\alpha\beta}u_{\alpha\beta}u_{jn}^2+I_j+J_{1j}+J_{2j}\right\}\nonumber\\
&&-\frac{1}{{\sigma}^3_1(B)}\sum_{\alpha,\beta=1}^n\sum_{i\in
B}F^{\alpha\beta}[{\sigma}_1(B)\widetilde{a}_{ii,\alpha}-\widetilde{a}_{ii}\sum_{j\in
B}\widetilde{a}_{jj,\alpha}][{\sigma}_1(B)\widetilde{a}_{ii,\beta}-\widetilde{a}_{ii}\sum_{j\in
B}\widetilde{a}_{jj,\beta}]\nonumber\\
&&-\frac{1}{{\sigma}_1(B)}\sum_{\alpha,\beta=1}^n\sum_{i\neq j,
i,j\in
B}F^{\alpha\beta}\widetilde{a}_{ij,\alpha}\widetilde{a}_{ij,\beta}+O(\mathcal
H_{\phi}).
\end{eqnarray}
\end{lemma}

\medskip

$I_j$, $J_{1j}$ and $J_{2j}$ in (\ref{3.37}) are crucial terms. Some
fine analysis of these terms are the key in our proof of the main
results.

\medskip

\noindent{\bf Proof of Lemma \ref{lem3.3}.} For any fixed point $x \in \O$, choose a
coordinate system as in (\ref{coord}) so that $|Du(x)|=  u_n(x) >0$
and the matrix $(\widetilde{a}_{ij}(x))$ is diagonal for $1 \le
i,j\le n-1$ and nonnegative.
From the definition of $p$,
\begin{eqnarray}\label{3.10}
-\frac{u_{jj}}{u_n}-\eta_0g=\widetilde{a}_{jj}= O(\mathcal H_{\phi}),
\forall j\in B; \quad p_{\alpha} =\sigma_{l}({G}) \sum_{j \in
B}\widetilde{a}_{jj,\alpha}+ O({\phi}).
\end{eqnarray}
By (\ref{3.10}),
\begin{eqnarray}\label{3.12a}
p_{\alpha\beta}&=& \sigma_l({G})[\sum_{j\in
B}\widetilde{a}_{jj,\alpha\beta}-2\sum_{i\in G,j\in
 B}\frac{\widetilde{a}_{ij,\alpha}\widetilde{a}_{ij,\beta}}{\widetilde{a}_{ii}}]
 +O(\mathcal H_{\phi})\nonumber\\
&=&\sigma_l({G})[\sum_{j\in
B}(a_{jj,\alpha\beta}-\eta_0g_{\alpha\beta})-2\sum_{i\in G,j\in
 B}\frac{\widetilde{a}_{ij,\alpha}\widetilde{a}_{ij,\beta}}{\widetilde{a}_{ii}}]
 +O(\mathcal H_{\phi}).
\end{eqnarray}

Since $u_k=0$ at $x$ for $k=1,\cdots, n-1$, from
(\ref{2.18}), and for each $j\in B$,
\begin{eqnarray}\label{3.19'}
u_n^{3}a_{jj,\alpha\beta}& =&-u_n^2u_{jj\alpha\beta
}-2u_n(u_{n\beta}u_{jj\alpha}+u_{n\a}u_{jj\b})+
2u_n(u_{j\alpha}u_{nj\b}+u_{j\b}u_{nj\a})\nonumber\\&&+2u_nu_{nj}u_{\a\b
j}+2u_{nj}(u_{n\a}u_{j\b}+u_{n\b}u_{j\a})-2u_{nn}u_{j\a}u_{j\b}
-(2u_{n\alpha}u_{n\beta}+\nonumber\\&&2u_nu_{\alpha\beta
n})u_{jj}-2\eta_0gu_{j\alpha}u_{j\beta}u_n-3\eta_0u_n^2(u_{n\alpha}g_{\beta}
+u_{n\b}g_{\a})\nonumber\\&&-\eta_0g(3u_n^2u_{n\alpha\beta}+6u_{n\alpha}
u_{n\beta}u_n+\sum_{i=1}^{n-1}u_{i\alpha}u_{i\beta}u_n)
 + O(\mathcal H_{\phi}).
\end{eqnarray}
From the definition of $a_{ij}$,
\begin{eqnarray}\label{3.21a}
u_n u_{ij\a}=-u_n^2 a_{ij,\a}+u_{nj}u_{i\a}
+u_{ni}u_{j\a}+u_{n\a}u_{ij}, \quad \forall \;i,j \le n-1,
\end{eqnarray}
and
\begin{eqnarray}\label{ta}\widetilde{a}_{ij,\a}=a_{ij\a}-\eta_0g_{\a}\delta_{ij}, \quad\forall\; i,j\in
B,\end{eqnarray}
\begin{eqnarray}\label{3.23}
&&\sum_{\alpha,\beta=1}^n F^{\alpha\beta}a_{jj,\alpha \beta}
\\&=& \sum_{\alpha,\beta=1}^n
\frac{F^{\alpha\beta}}{u_n^3}[-u_n^2u_{\alpha\beta jj}-4u_{n\a}u_{nj}u_{j\b}
+4u_nu_{j\a}u_{nj\b} +2u_nu_{nj}u_{\a\b
j}-2u_{nn}u_{j\a}u_{j\b}\nonumber\\
&&-\eta_0g(2u_{j\a}u_{j\b}u_n+u_{n\a\b}u_n^2+\sum_{i=1}^{n-1}u_{i\a}u_{i\b}u_n)
-2\eta_0u_{n\a}g_{\b}u_n^2] +O(\mathcal H_{\phi}).\nonumber
\end{eqnarray}
Break summation as $\sum_{\alpha=1}^n F^{\alpha
n}u_{n\a}=(\sum_{\alpha,\beta=1}^n -\sum_{\b
=1}^{n-1}\sum_{\a=1}^n)F^{\alpha\beta}u_{\a\b}$, $\forall j\in B$,
\begin{eqnarray*}\label{3.24}
\sum_{\alpha,\beta=1}^n
F^{\alpha\beta}u_{n\a}u_{j\b}=u_{nj}(\sum_{\alpha,\beta=1}^n
- \sum_{\b
=1}^{n-1}\sum_{\a=1}^n)F^{\alpha\beta}u_{\a\b}+\sum_{\a=1}^nF^{j\a}u_{n\a}u_{jj},
\end{eqnarray*}
\begin{eqnarray*}\label{3.25}
\sum_{\alpha,\beta=1}^n
F^{\alpha\beta}u_{j\a}u_{nj\b}
=u_{nj}(\sum_{\alpha,\beta=1}^n
-\sum_{\a=1}^{n}\sum_{\b=1}^{n-1})F^{\alpha\beta}u_{\a\b
j}+\sum_{\a=1}^nF^{j\a}u_{j\a n}u_{jj},
\end{eqnarray*}
and
by (\ref{3.21a}), for $j\in B$,
\begin{eqnarray*}\label{3.28}
&& u_n\sum_{\a=1}^n\sum_{\b=1}^{n-1}F^{\alpha\beta}u_{\a\b
j}=u_n\sum_{\a=1}^{n}\bigg(\sum_{i \in
B}F^{\alpha i}u_{ij\a}+\sum_{i \in G}F^{\alpha
i}u_{ij\a}\bigg)\nonumber\\
&=&\sum_{\a=1}^{n}\sum_{i \in B}F^{\alpha
i}(u_{i\a}u_{jn} + u_{j\a}u_{in}-\eta_0g_{\a}u_n^2\delta_{ij}+u_{n\a}u_{ij}) \nonumber\\
&&+\sum_{\a=1}^{n} \sum_{i \in G}F^{\alpha i}
(-u_n^2\widetilde{a}_{ij,\a}+u_{i\a}u_{jn} + u_{j\a}u_{in})
+O(\mathcal H_{\phi})\nonumber\\
&=&-u_n^2\sum_{\a=1}^{n} \sum_{i \in G}F^{\alpha i}
\widetilde{a}_{ij,\a}+u_{nj}\sum_{i\in G}F^{ii}u_{ii}
+2u_{nj}(\sum_{i=1}^{n-1}F^{ni}u_{ni})
+u_{nj}\sum_{i\in B}F^{ii}u_{ii}\nonumber\\
&&+u_{jj}\sum_{i=1}^{n-1}F^{ij}u_{ni}
+u_{jj}\sum_{i=1}^{n}F^{ij}u_{ni}-\eta_0\sum_{\a=1}^nF^{j\a}g_{\a}u_n^2
+O(\mathcal H_{\phi}),
\end{eqnarray*}
and
\begin{eqnarray*}
&&\sum_{\alpha,\beta=1}^n
F^{\alpha\beta}u_{j\a}u_{j\b}= F^{nn}u_{nj}^2+2F^{jn}u_{jn}u_{jj}+F^{jj}u_{jj}^2\\
&=&u_{nj}^2(\sum_{\alpha,\beta=1}^n
F^{\alpha\beta}u_{\a\b}-2\sum_{\a=1}^{n-1}F^{\a
n}u_{n\a}-\sum_{\alpha,\beta=1}^{n-1}
F^{\alpha\beta}u_{\a\b})+2F^{jn}u_{jn}u_{jj}+F^{jj}u_{jj}^2.
\end{eqnarray*}
Put above to (\ref{3.23}),
\begin{eqnarray}\label{3.29}
u_n^3\sum_{\alpha,\beta=1}^nF^{\alpha\beta} a_{jj,\alpha\beta}
&=&-u_n^2\sum_{\alpha,\beta=1}^nF^{\alpha\beta}(u_{\alpha\beta jj}
+6u_{\a\b j}u_{nj}u_n-6u_{\a\b}u_{nj}^2)\nonumber\\
&& +4u_n^2u_{nj}\sum_{\a=1}^{n} \sum_{i \in G}F^{\alpha
i}\widetilde{a}_{ij,\a}
+2u_{nj}^2\sum_{i\in G}F^{ii}u_{ii}+2u_{nj}^2\sum_{i\in B}F^{ii}u_{ii}\nonumber\\
&&-12\sum_{\a=1}^nF^{j\a}u_{n\a}u_{jn}u_{jj}+4\sum_{\a=1}^nF^{j\a}u_{jn\a}u_{n}u_{jj}
-2u_{nn}F^{jj}u_{jj}^2\nonumber\\
&&-\eta_0g\sum_{\alpha,\beta=1}^nF^{\alpha\beta}(2u_{j\a}u_{j\b}u_n+u_{n\a\b}u_n^2+
\sum_{i=1}^{n-1}u_{i\a}u_{i\b}u_n)\nonumber\\
&&-2\eta_0\sum_{\alpha,\beta=1}^n F^{\alpha\beta}u_{n\a}g_{\b}u_n
+4\eta_0\sum_{\a=1}^nF^{j\a}g_{\a}u_{jn}u_n^2+O(\mathcal H_{\phi}).
\end{eqnarray}

Since $a_{ij,\a}=\widetilde{a}_{ij,\a}$ for $i \neq j$,
\begin{eqnarray}\label{3.29-0}
\sum_{\alpha,\beta=1}^n
\dfrac{F^{\alpha\beta}{\widetilde{a}_{ij,\alpha}\widetilde{a}_{ij,\beta}}}{\widetilde{a}_{ii}}
=\sum_{\alpha,\beta=1}^n
\dfrac{F^{\alpha\beta}{{a}_{ij,\alpha}{a}_{ij,\beta}}}{{a}_{ii}}\frac{a_{ii}}{\widetilde{a}_{ii}}
=\sum_{\alpha,\beta=1}^n
\dfrac{F^{\alpha\beta}{{a}_{ij,\alpha}{a}_{ij,\beta}}}{{a}_{ii}}(1+\frac{\eta_0g}{\widetilde{a}_{ii}}).
\end{eqnarray}
(\ref{3.12a}), (\ref{3.29}) and \eqref{3.29-0} yield that, for each $j\in B$,
\begin{eqnarray}\label{3.30}
F^{\a\b}p_{\a\b} &=&F^{\a\b}\sigma_l({G})[\sum_{j\in
B}(a_{jj,\alpha\beta}-\eta_0g_{\alpha\beta})-2\sum_{i\in G,j\in
 B}\frac{\widetilde{a}_{ij,\alpha}\widetilde{a}_{ij,\beta}}{\widetilde{a}_{ii}}]
 +O(\mathcal H_{\phi})\nonumber
\\&=&
u_n^{-3}\sum_{j \in B}\sigma_l(G)\Big[-\sum_{\alpha,\beta=1}^n
F^{\alpha\beta}u_n^2u_{\alpha\beta jj}+6 u_n
\sum_{\alpha,\beta=1}^nF^{\alpha\beta}u_{jn}u_{\alpha\beta
j}\nonumber\\
&&-6u_{jn}^2 \sum_{\alpha,\beta=1}^nF^{\alpha\beta}u_{\alpha \beta}+
I_j+J_{1j}\Big]  +O(\mathcal H_{\phi}).
\end{eqnarray}
where $I_j, J_{1j}$ as in \eqref{I-j} and \eqref{Jp}.

For each $j\in B$, differentiating equation (\ref{1.1}) in $e_j$
direction at $x$,
\begin{equation}\label{4.10}
\sum_{\alpha,\beta=1}^nF^{\alpha\beta}u_{\alpha\beta
j}+F^{u_n}u_{jn}+F^{u_j}u_{jj}+F^{x_j}=0,
\end{equation}
\begin{eqnarray}\label{4.11}
-\sum_{\alpha,\beta=1}^nF^{\alpha\beta}u_{\alpha\beta jj}
&=&\sum_{\alpha,\beta,\g,\e=1}^nF^{\alpha\beta, \g\e}u_{\alpha\beta
j}u_{\g\e j} +2\sum_{\alpha,\beta,l=1}^nF^{\alpha\beta,
u_l}u_{\alpha\beta j}u_{lj}\\
&&+2\sum_{\alpha,\beta=1}^nF^{\alpha\beta, x_j}u_{\alpha\beta
j}+\sum_{l,s=1}^nF^{u_l, u_s}u_{lj}u_{sj}+F^{u}u_{jj}\nonumber\\
&&+\sum_{l=1}^nF^{u_l,
x_j}u_{lj}+F^{x_j,x_j}+\sum_{l=1}^nF^{u_l}u_{ljj}.\nonumber
\end{eqnarray}
It follows from (\ref{3.21a}) that, at $x$
\begin{eqnarray}\label{4.12}
-\sum_{\alpha,\beta=1}^nF^{\alpha\beta}u_{\alpha\beta
jj}&=&\sum_{\alpha,\beta,\g,\e=1}^nF^{\alpha\beta,
\g\e}u_{\alpha\beta j}u_{\g\e
j}+2\sum_{\alpha,\beta=1}^n[F^{\alpha\beta,
u_n}u_{j\alpha\beta}u_{nj}+F^{\alpha\beta,
x_j}u_{\alpha\beta j}]\nonumber\\
&&+F^{u_n, u_n}u^2_{jn}+2F^{u_n,
x_j}u_{jn}+F^{x_j,x_j}+2\frac{F^{u_n}}{u_n}u^2_{jn}+J_{2j}
+O(\mathcal H_{\phi}),
\end{eqnarray}
where $J_{2j}$ is defined in \eqref{J2}.

Since $u_{\alpha\beta jj}=u_{ jj\alpha\beta}$, from
(\ref{3.30}) and (\ref{4.12}),
\begin{eqnarray}\label{pab}
F^{\a\b}p_{\a\b} &=&\sum_{j\in B}u_n^{-3}\sigma_l(G)\Big[
 \Big(\sum_{\alpha,\beta,\g,\e=1}^n
F^{\alpha\beta, \g\e}u_{\alpha\beta j}u_{\g\e j}
+2\sum_{\alpha,\beta=1}^nF^{\alpha\beta,
u_n}u_{j\alpha\beta}u_{jn}\nonumber\\
&&+2\sum_{\alpha,\beta=1}^nF^{\alpha\beta,
x_j}u_{j\alpha\beta}+F^{u_n,u_n}u_{jn}^2+2F^{u_n,x_j}u_{jn}+F^{x_j,
x_j}+2\frac{F^{u_n}}{u_n}u_{jn}^2\Big)u_n^2\nonumber\\
&&+6\sum_{\alpha,\beta=1}^nF^{\alpha\beta}u_{j\alpha\beta}u_{jn}u_n
-6\sum_{\alpha,\beta=1}^nF^{\alpha\beta}u_{\alpha\beta}u_{jn}^2+I_j+J_{1j}+J_{2j}\Big]
+O(\mathcal H_{\phi}).
\end{eqnarray}

The fact $q\in C^{1,1}(\O)$ follows
Corollary 2.2 in \cite{BG08}. Also by Lemma 2.4 in
\cite{BG08},
\begin{equation}\label{3.36}
q_\alpha=\frac{\partial q}{\partial x_\alpha}=\sum_{j \in
B}\frac{{\sigma}^2_1(B|j)-{\sigma}_2(B|j)}{{\sigma}^2_1(B)}\widetilde{a}_{jj,\alpha}
+O(\mathcal H_{\phi}),
\end{equation}
and
\begin{eqnarray}\label{q-eu}
q_{\alpha \beta}&=&\sum_{j\in
B}\frac{{\sigma}^2_1(B|j)-{\sigma}_2(B|j)}{{\sigma}^2_1(B)}\Big[\widetilde{a}_{jj,\alpha
\beta}-2\sum_{i\in
G}\frac{\widetilde{a}_{ij,\alpha}\widetilde{a}_{ij,\beta}}{\widetilde{a}_{ii}}\Big]
\nonumber\\
&&-\frac{1}{{\sigma}^3_1(B)}\sum_{i\in
B}\Big[{\sigma}_1(B)\widetilde{a}_{ii,\alpha}-\widetilde{a}_{ii}\sum_{j\in
B}\widetilde{a}_{jj,\alpha}\Big]\Big[{\sigma}_1(B)\widetilde{a}_{ii,\beta}-\widetilde{a}_{ii}\sum_{j\in
B}\widetilde{a}_{jj,\beta}\Big]\nonumber\\
&&-\frac{1}{{\sigma}_1(B)}\sum_{i\neq j\in
B}\widetilde{a}_{ij,\alpha}\widetilde{a}_{ij,\beta}+O(\mathcal
H_{\phi}).
\end{eqnarray}
Following the same computations as for $p$, we get
\begin{eqnarray}
&&\sum_{\alpha,\beta=1}^{n}F^{\alpha\beta}q_{\alpha\beta}\nonumber\\
&=&\sum_{j\in
B}\frac{{\sigma}^2_1(B|j)-{\sigma}_2(B|j)}{{\sigma}^2_1(B)u_n^{3}}\Big[
 \Big(\sum_{\alpha,\beta,\g,\e=1}^n
F^{\alpha\beta, \g\e}u_{\alpha\beta j}u_{\g\e
j}+2\sum_{\alpha,\beta=1}^nF^{\alpha\beta,
u_n}u_{j\alpha\beta}u_{jn}\nonumber\\
&&+2\sum_{\alpha,\beta=1}^nF^{\alpha\beta,
x_j}u_{j\alpha\beta}+F^{u_n,u_n}u_{jn}^2+2F^{u_n,x_j}u_{jn}+F^{x_j,
x_j}+2\frac{F^{u_n}}{u_n}u_{jn}^2\Big)u_n^2\nonumber\\
&&+6\sum_{\alpha,\beta=1}^nF^{\alpha\beta}u_{j\alpha\beta}u_{jn}u_n
-6\sum_{\alpha,\beta=1}^nF^{\alpha\beta}u_{\alpha\beta}u_{jn}^2+I_j+J_{1j}+J_{2j}\Big]\nonumber\\
&&-\frac{1}{{\sigma}^3_1(B)}\sum_{\alpha,\beta=1}^n\sum_{i\in
B}F^{\alpha\beta}[{\sigma}_1(B)\widetilde{a}_{ii,\alpha}-\widetilde{a}_{ii}\sum_{j\in
B}\widetilde{a}_{jj,\alpha}][{\sigma}_1(B)\widetilde{a}_{ii,\beta}-\widetilde{a}_{ii}\sum_{j\in
B}\widetilde{a}_{jj,\beta}]\nonumber\\
&& -\frac{1}{{\sigma}_1(B)}\sum_{\alpha,\beta=1}^n\sum_{i\neq j\in
B}F^{\alpha\beta}\widetilde{a}_{ij,\alpha}\widetilde{a}_{ij,\beta}+O(\mathcal
H_{\phi}).
\end{eqnarray}
The proof of the Lemma is complete. \qed

\section{Proof of Theorems}

We want to show $\widetilde{a}$ defined in (\ref{a}) is of constant rank. Theorem \ref{th1.1} corresponds
to the case $\eta_0=0$. As for Theorem \ref{th1.2}, set
\begin{equation}\label{stress} \Omega_{\varpi}=\{x \in \Omega | 0< \k_s(x)< \frac{\lambda}{100\varpi}.\}\end{equation}
We indicate how the constant rank theorem for $\widetilde{a}$ would imply Theorem \ref{th1.2}. If $\min\{\k^0, \k^1\}=0$, the strict
convexity of level surfaces $\Sigma^c$ for $c\in (0,1)$ in Theorem
\ref{th1.2} follows from Theorem \ref{th1.1}. We may assume $\min\{\k^0, \k^1\}>0$.  By Theorem
\ref{th1.2}, $a$ is strictly positive definite in $\bar \Omega$. That
is, $\k_s(x)>0, \forall x\in \bar \Omega$. By the continuity, $\k_s(x)$ has a positive lower bound
(which we want to estimate).
Increasing $\eta_0$ from $0$ to the level that $\widetilde{a}$ is nonnegative definite through out $\bar \Omega_{\varpi}$ but degenerate at some points $x_0$. (\ref{k-est}) follows easily if the degeneracy happens
outside $\Omega_{\varpi}$ or on the boundary.  If the degeneracy
happens at an interior point $x_0$ of $\Omega_{\varpi}$, the goal is to show that $\widetilde{a}$ is degenerate
through out the connected component $U$ of $\Omega_{\varpi}$ containing $x_0$ with the same rank. If this is true, $\k_s(x)=constant, \forall x\in \Sigma^c\bigcap U$. Theorem \ref{th1.2} would follow. Furthermore,
the closure of $\Sigma^{u(x)}\bigcap U$ can not intersect the set $\{z\in \Omega |\k_s(z)=\frac{\lambda}{100\varpi}\}$ for any $x\in U$, since by definition $\k_s(x)<\frac{\lambda}{100\varpi}$.
This would imply that if the degeneracy of $\widetilde{a}$
happens at an interior point of $\Omega_{\varpi}$, the connect component $U$ of $\Omega_{\varpi}$ containing
that point is exactly the set $\bigcup_{c_0<c<c_1} \Sigma^c$ for some $0\le c_0<c_1\le 1$, and $\k_s(x)\equiv constant, \forall x\in \Sigma^c$. Therefore, $\Sigma^c$ is a round sphere for each $c_0\le c\le c_1$.
The task now is to prove the rank of $\widetilde{a}$ is constant. That is a consequence of the following proposition.

\begin{proposition}\label{thm4.1}
Suppose $u\in C^{3,1}$ is a quasiconcave solution of equation (\ref{1.1}) and $F$ satisfies assumptions in Theorem
\ref{th1.1}. If the second
fundamental form of $\Sigma^c$ of solution $u$ attains minimum rank
$l$ at certain point $x_0\in \Omega$, then there exist a
neighborhood $\O$ of $x_0$ and a positive constant $C$ independent
of $\phi$ (defined in (\ref{3.4}) ), such that
\begin{equation}\label{4.4}
\sum_{\alpha,\beta=1}^{n}F^{\alpha\beta}\phi_{\alpha\beta}(x)\leq
C(\phi(x)+|\nabla \phi(x)|),~~\forall ~x\in\O.
\end{equation}
If in addition $F$ satisfies the uniform ellipticity condition
(\ref{1.2-0}) in $\O\subset \Omega_{\varpi}$, then there is $A$
depending only on $||F||_{C^{2}}, n, \lambda, d_0, \|u\|_{C^4(\bar
\Omega)}$, such that if $\widetilde a(x)$ defined in (\ref{a})
attains minimum rank $l$ at certain point $x_0\in \O$, then
inequality (\ref{4.4}) is true for all $x\in \O$.
\end{proposition}

\noindent{\bf Proof of Proposition~\ref{thm4.1}.}  Suppose
the minimum rank $l$ of $\widetilde a$ is attained at an interior point $x_0$, and we may assume $l\le n-2$. Let $\O$ be a small neighborhood of $x_0$. Lemma \ref{lem3.3} and (\ref{3.4})
implies $\phi\in C^{1,1}(\O),$ $ \phi(x)\geq 0, \;\phi(x_0)=0 $. For
$\epsilon>0$ sufficient small, let $\phi_{\epsilon}$ defined as in
(\ref{3.5}), we want to establish differential inequality
(\ref{4.4}) for $\phi_\varepsilon$  with constant $C$ independent of
$\varepsilon$ in $\O$.  For each fixed $x\in \O$, choose a local
coordinate frame $e_1,\cdots, e_{n-1}, e_n$  so (\ref{coord}) and
(\ref{G-B}) are satisfied. We will omit the
subindex $\varepsilon$ with the understanding that all the estimates
are independent of $\varepsilon$.
From Lemma~\ref{lem3.3},

\begin{eqnarray}\label{4.13}
&&F^{\alpha\beta}\phi_{\alpha\beta}\nonumber\\&=& \sum_{j\in
B}u_n^{-3}\left[\sigma_l(G)+\frac{{\sigma}^2_1(B|j)-{\sigma}_2(B|j)}{{\sigma}^2_1(B)}\right]\left\{
 S_j+I_j+J_{1j}+J_{2j}\right\}\nonumber\\
&&-\frac{1}{{\sigma}^3_1(B)}\sum_{\alpha,\beta=1}^n\sum_{i\in
B}F^{\alpha\beta}[{\sigma}_1(B)\widetilde{a}_{ii,\alpha}-\widetilde{a}_{ii}\sum_{j\in
B}\widetilde{a}_{jj,\alpha}][{\sigma}_1(B)\widetilde{a}_{ii,\beta}-\widetilde{a}_{ii}\sum_{j\in
B}\widetilde{a}_{jj,\beta}]\nonumber\\
&&-\frac{1}{{\sigma}_1(B)}\sum_{\alpha,\beta=1}^n\sum_{i\neq j,
i,j\in
B}F^{\alpha\beta}\widetilde{a}_{ij,\alpha}\widetilde{a}_{ij,\beta}+O(\mathcal
H_{\phi}),
\end{eqnarray}
where $I_j, J_{1j}, J_{2j}$ are defined in (\ref{I-j}), (\ref{Jp}), (\ref{J2}) respectively, and,
\begin{eqnarray}\label{4.14}
S_j= && \Big[\sum_{\alpha,\beta,\g,\e=1}^n F^{\alpha\beta,
\g\e}u_{j\alpha\beta}u_{\g\e
j}+2\sum_{\alpha,\beta=1}^nF^{\alpha\beta,
u_n}u_{j\alpha\beta}u_{jn}+2\sum_{\alpha,\beta=1}^nF^{\alpha\beta,
x_j}u_{j\alpha\beta}\nonumber\\
&+&F^{u_n,u_n}u_{jn}^2+2F^{u_n,x_j}u_{jn}+F^{x_j,
x_j}+2\frac{F^{u_n}}{u_n}u_{jn}^2\Big]u_n^2\nonumber\\
&+&6\sum_{\alpha,\beta=1}^nF^{\alpha\beta}u_{j\alpha\beta}u_{jn}u_n
-6\sum_{\alpha,\beta=1}^nF^{\alpha\beta}u_{\alpha\beta}u_{jn}^2
\end{eqnarray}

In the coordinate system (\ref{coord}),
\begin{equation}\label{coord-t}t=u_n^{-1}, \quad
D^2u(x)=t^{-1}\mathbf   {\widetilde A}, \quad A_{ij}=tu_{ij}=\frac{u_{ij}}{u_n},
\quad \t=(0,..,0,1).\end{equation}

For each $j \in B$, set
\begin{eqnarray}\label{V}\begin{array}{rcl}
&&\widetilde{X}_{\a\b}=u_{\a\b j}u_n ;\quad\forall \;\a,\b \in G\cup
\{n\}\quad \text{with}\quad(\alpha, \beta)\neq (n,n),\\
&&\widetilde{X}_{\a\b}=2u_{\a\b}u_{jn}; \quad \forall\; \a\in B\;
\text{or}\; \forall\; \b \in B,
\\&&\widetilde{X}_{nn}
=u_{nn j}u_n+\dfrac{1}{F^{nn}}[2\sum_{\substack{\a \in B\\ \b \in
G\cup\{n\}}}F^{\a\b}u_{\a\b j}u_n+\sum_{\a,\b \in B}F^{\a\b}u_{\a\b j}u_n\\
&&\quad\quad\quad-2\sum_{\a,\b \in
B}F^{\a\b}u_{\a\b}u_{jn}-4\sum_{\a \in B}F^{\a n}u_{\a n}u_{jn}
+F^{u_j}u_{jj}u_n];\\
&&\widetilde{Y}=u_{jn}u_n; \quad Z_k=\d_{kj}u_n.
\end{array}
\end{eqnarray}
For such $\widetilde V=\big((\widetilde{X}_{\a\b}), \widetilde Y,
(Z_k)\big)$, by \eqref{4.10},
\begin{eqnarray*}\label{Vb}
<V,\nabla_{(B,x)}F>&=&{F^{\a\b}}\widetilde{X}_{\a\b}+F^{u_n}\widetilde{Y}+F^{x_k}Z_k\\
&=&F^{\a\b}u_{\a\b j}u_n
+F^{u_j}u_{jj}u_n+F^{u_n}u_{jn}u_n+F_{x_j}u_n=0.
\end{eqnarray*} We need the following lemma.

\begin{lemma}\label{Claim} Under the coordinate system (\ref{coord}) at
$x$ with $\widetilde{V}$ as in \eqref{V}, then
\[u_n^3
I=I_j+O(u_{jj}),\] where $I_j$ defined in (\ref{I-j}) and $I$ defined
in (\ref{sdQZ-0}).
\end{lemma}

\noindent {\bf Proof of Lemma \ref{Claim}.} Since
 $u_{ii}=-{a}_{ii}u_n$ for $i \in G$ and
$(u_{ij}),i,j=1,..,n-1$ is diagonal at $x$. By (\ref{I-j}) and \eqref{3.21a}, for
each $j \in B$,
 \begin{eqnarray}\label{Ij}
I_j&=&\sum_{i\in G}[-2u_n^3\sum_{\a,\b\notin B }
\dfrac{F^{\alpha\beta}{a}_{ij,\alpha}{a}_{ij,\beta}}{{a}_{ii}}
+4u_n^{2}u_{nj}\sum_{\a \notin B } F^{\alpha
i}{a}_{ij,\a}+2u_{nj}^2F^{ii}{u}_{ii}]\nonumber\\
&=&\sum_{i\in G}\sum_{\a,\b\notin B }[2u_n^2
\dfrac{F^{\alpha\beta}{u}_{ij,\alpha}{u}_{ij,\beta}}{{u}_{ii}}
-8u_nu_{nj}\frac{ F^{\alpha
\b}u_{ij,\a}u_{i\b}}{u_{ii}}+8u_{nj}^2\frac{F^{\a\b}u_{i\a}u_{i\b}}{u_{ii}}]+O(u_{jj})\nonumber\\
&=&2\sum_{i\in G}\sum_{\a,\b\notin B
}\frac{F^{\a\b}(u_{ij\a}u_n-2u_{i\a}u_{jn})(u_{ij\b}u_n-2u_{i\b}u_{jn})}{u_{ii}}+O(u_{jj}).
\end{eqnarray}

By (\ref{coord-t}), \eqref{newX1}, (\ref{Y-d2}) and \eqref{V}, \[Y_{i\a}=0, \forall\; i \in B; \quad
Y_{i\a}=-(u_{ij\a}u_n^{-1}-2u_{i\a}u_{jn}u_n^{-2}), \forall \;i \in
G , \a \notin B.\] Therefore
\begin{eqnarray}\label{I}
u_n^3 I&=&2\sum_{i \in G}\sum_{\a,\b\notin B
}\frac{F^{\a\b}}{u_{ii}}Y_{i\a}Y_{i\b}u_n^4=I_j+O(u_{jj}).
\end{eqnarray}  \qed

\medskip

We now finish the proof of the proposition.

In case of \;{\bf Theorem~\ref{th1.1}}, since $\eta_0=0$,
$u_{jj}=O(\phi), \forall \;j \in B$, $a=\widetilde{a}$. \eqref{4.13}
gives
\begin{eqnarray}\label{phi}
&&F^{\alpha\beta}\phi_{\alpha\beta}\nonumber\\&=& \sum_{j\in
B}u_n^{-3}\left[\sigma_l(G)+\frac{{\sigma}^2_1(B|j)-{\sigma}_2(B|j)}{{\sigma}^2_1(B)}\right]
 (S_j+I_j)\nonumber\\
&&-\frac{1}{{\sigma}^3_1(B)}\sum_{\alpha,\beta=1}^n\sum_{i\in
B}F^{\alpha\beta}[{\sigma}_1(B){a}_{ii,\alpha}-{a}_{ii}\sum_{j\in
B}{a}_{jj,\alpha}][{\sigma}_1(B){a}_{ii,\beta}-{a}_{ii}\sum_{j\in
B}{a}_{jj,\beta}]\nonumber\\
&&-\frac{1}{{\sigma}_1(B)}\sum_{\alpha,\beta=1}^n\sum_{i\neq j,
i,j\in B}F^{\alpha\beta}{a}_{ij,\alpha}{a}_{ij,\beta}+O(\mathcal
H_{\phi}).
\end{eqnarray}

In fact for each $j \in B$, by \eqref{V} and
\eqref{3.21a},
\begin{eqnarray}\label{n0}\begin{array}{rcl}
&&\widetilde{X}_{\a\b}=u_{\a \b j}u_{n}+O(\mathcal H_{\phi}); \quad
\forall\; \a\in B\; \text{or}\; \forall\; \b \in B,\\
&&\widetilde{X}_{\a\b}=u_{\a \b j}u_{n}; \quad
\forall\; \a, \b\in G\cup\{n\}, \, (\a,\b)\neq (n,n),
\\&&\widetilde{X}_{nn}
=u_{nn j}u_n+O(\mathcal H_{\phi}).
\end{array}
\end{eqnarray}

Lemma \ref{Claim}, condition \eqref{1.4} and
Lemma~\ref{lem3.5-new} imply
\begin{equation}\label{Sj} S_j+I_j=H(\widetilde V,\widetilde V)+O(\mathcal H_{\phi})\le O(\mathcal H_{\phi}) .
\end{equation}

By the Newton-MacLaurine inequality, $C \ge
\sigma_l(G)+\frac{{\sigma}^2_1(B|j)-{\sigma}_2(B|j)}{{\sigma}^2_1(B)}
\ge 0$. Condition (\ref{1.2}) implies that there is $\delta>0$
\begin{equation}\label{4.15}
 (F^{\alpha\beta}) \ge \delta I, \; \forall \;x \in \O.
\end{equation}
Combining (\ref{phi}), (\ref{Sj}) and (\ref{4.15})
\begin{eqnarray}\label{4.17}
F^{\alpha\beta}\phi_{\alpha\beta}\le && C(\phi+\sum_{i,j\in
B}|\nabla a_{ij}|)-\frac{\delta}{\sigma_1(B)}\sum_{
\alpha=1}^n\sum_{i\neq j\in B}a^2_{ij\alpha}\\\nonumber
&&-\frac{\delta}{{\sigma}^3_1(B)}\sum_{\alpha=1}^n\sum_{i\in
B}({\sigma}_1(B)a_{ii,\alpha}-a_{ii}\sum_{j\in B}a_{jj,\alpha})^2.
\end{eqnarray}
Finally, by Lemma~3.3 in \cite{BG08}, the term $\sum_{i,j \in B}
|\nabla a_{ij}|$  can be controlled by the rest terms on the right
hand side in (\ref{4.17}) and $\phi + |\nabla \phi|$. In conclusion,
there exist positive constant $C$ independent of
$\varepsilon$, such that
\begin{eqnarray}\label{4.19}
\sum_{\alpha,\beta}
F^{\alpha\beta}\phi_{\alpha\beta}&\le&C(\phi+|\nabla \phi|).
\end{eqnarray}
{\bf Proposition \ref{thm4.1}} is verified under the assumptions of \;{\bf
Theorem~\ref{th1.1}}. Therefore, Theorem \ref{th1.1} is proved.

\medskip

The rest proof is to deal with the case $\eta_0>0$. In this case, since $\k_s(x)>0$ for all $x\in \Omega$
by Theorem \ref{th1.1}. Since $\O \subset \Omega_{\varpi}$,
\begin{equation}\label{stress1} 0< \k_s(x)< \frac{\lambda}{100\varpi}, \forall x\in \O.\end{equation}

For each $j \in B$, it follows from \eqref{V}, \eqref{3.21a}, \eqref{ta}, and the fact
$u_k=0$ for all $k\le n-1$,
\begin{eqnarray}\label{n00}\begin{array}{rcl}
\widetilde{X}_{\a\b}
&=&u_{\a \b j}u_{n}+\widetilde{a}_{\a
j\b}u_n^2+\eta_0g_{n}u_n^2\delta_{\b n}\delta_{\a j}+O(u_{jj}), \quad \forall \a\in
B;\\
\widetilde{X}_{\a\b}&=&u_{\a \b j}u_{n}; \quad
\forall\; \a, \b\in G\cup\{n\}, \, (\a,\b)\neq (n,n),
\\
\widetilde{X}_{nn}
&=&u_{nn j}u_n-\frac{2\eta_0g_nu_n^2F^{nj}}{F^{nn}}+ O(\sum_{i,j\in
B}|\nabla \widetilde{a}_{ij}|)+O(u_{jj}).
\end{array}
\end{eqnarray}

In view of (\ref{n00}), \begin{eqnarray}\label{n000}\begin{array}{rcl}
u_{\a \b j}u_{n}
&=&\widetilde{X}_{\a\b}-\widetilde{a}_{\a
j\b}u_n^2-\eta_0g_{n}u_n^2\delta_{\b n}\delta_{\a j}+O(u_{jj}), \quad \forall \a\in
B;\\
u_{\a \b j}u_{n}&=&\widetilde{X}_{\a\b}; \quad
\forall\; \a, \b\in G\cup\{n\}, \, (\a,\b)\neq (n,n),
\\
u_{nn j}u_n
&=&\widetilde{X}_{nn}+\frac{2\eta_0g_nu_n^2F^{nj}}{F^{nn}}+ O(\sum_{i,j\in
B}|\nabla \widetilde{a}_{ij}|)+O(u_{jj}).
\end{array}
\end{eqnarray}
Notice that $u_{jj}=-\eta_0gu_n+O(\phi), \; \forall \; j \in B$. Substitute $u_{\a \b j}u_{n}$
by formula (\ref{n000}) in $S_{j}$
defined in (\ref{4.14}). We need to track the terms with factor $\eta_0^2g_n^2$. They are coming from $\sum_{\a, \b, \g, \e
\in\{j,n\}} F^{\alpha\beta, \g\e}u_{j\alpha\beta}u_{\g\e j}u_n^2$ only.
In turn, the coefficient in front of $\eta_0^2g_n^2$ can be controlled by,
say $50\varpi\eta_0^2 u_n^3$, where $\varpi$ is
defined in (\ref{varpi}). By Lemma \ref{Claim}, condition
\eqref{1.4}, Lemma~\ref{lem3.5-new} and the assumptions of Theorem
\ref{th1.2}, there exist constants $C_0', C_0''$ depending only on
$\lambda, \|F\|_{C^2}, d_0, \|u\|_{C^3}$ such that
\begin{eqnarray}\label{new-000}
S_j+I_j&\le& H(\widetilde V,\widetilde V)+50\varpi\eta_0^2 g_n^2u_n^3 +C_0'\eta_0|g_{n}|+C_0''\eta_0g+O(\mathcal H_{\phi})\nonumber\\
&\le& 50\varpi\eta_0^2 g_n^2u_n^3+C_0'\eta_0|g_{n}|
+C_0''\eta_0g+O(\mathcal H_{\phi})\nonumber\\
&=& (50\varpi \eta_0 g) A^2\eta_0 g
u_n^5+C_0'A\eta_0gu_{n}
+C_0''\eta_0g+O(\mathcal H_{\phi}).\end{eqnarray}  Since $F^{nn}\ge \lambda$, by \eqref{Jp} and \eqref{J2},
\begin{eqnarray}\label{newJ} J_{1j}+J_{2j}+O(u_{jj}) &\le& -\eta_0
F^{nn}g_{nn}u_n^3+\eta_0C_1'|g_{n}|+\eta_0C_2'g+O(\phi)\nonumber\\
&\le&-\eta_0g \lambda A^2u_n^5+\eta_0g
AC_1''u_{n}+\eta_0C_2'g+O(\phi), \end{eqnarray} where $C_1', C_1''
C_2'$ are positive constants depending only on $n, d_0,
\;||u||_{C^{3}(\Omega)}, \|F\|_{C^2}$. Note that $\eta_0g\le
\k_s(x)<\frac{\lambda}{100\varpi}$, it follows from (\ref{stress1}) that
\begin{eqnarray}\label{h}S_j+I_j+J_{1j}+J_{2j}+O(u_{jj}) \le
\frac{\eta_0g}2(-A^2\lambda u_n^5+C_3'A+C_3'')+O(\mathcal H_{\phi}).
\end{eqnarray} where  $C_3',C_3''$ are
positive constants depending only on $n, d_0, \;||u||_{C^{3}(\Omega)}, \|F\|_{C^2}$.
Since $u_n\ge d_0>0$, we may choose $A$ large enough in \eqref{h} depending only on $n, d_0, \;||u||_{C^{3}(\Omega)}, \|F\|_{C^2}, \lambda$ such that
\begin{equation}\label{sign}S_j+I_j+J_{1j}+J_{2j}+O(u_{jj}) \le
O(\mathcal H_{\phi}).\end{equation}
By (\ref{4.13}) and (\ref{sign}), \begin{eqnarray}\label{4.17b}
F^{\alpha\beta}\phi_{\alpha\beta}\le && C(\phi+\sum_{i,j\in
B}|\nabla \tilde a_{ij}|)-\frac{\delta}{\sigma_1(B)}\sum_{
\alpha=1}^n\sum_{i\neq j\in B}\tilde a^2_{ij\alpha}\\\nonumber
&&-\frac{\delta}{{\sigma}^3_1(B)}\sum_{\alpha=1}^n\sum_{i\in
B}({\sigma}_1(B)\tilde a_{ii,\alpha}-\tilde a_{ii}\sum_{j\in B}\tilde a_{jj,\alpha})^2.
\end{eqnarray} As in the case of $\eta_0=0$, the same argument yields (\ref{4.19}) for $\tilde a$. Thus {\bf Proposition \ref{thm4.1}} is validated
under the assumptions in {\bf Theorem~\ref{th1.2}}. \qed

\medskip

We remark that Theorem \ref{th1.2} covers all quasilinear equations satisfying structural conditions (\ref{1.4})-(\ref{1.2-0}). Therefore, it covers the quasilinear equations treated in \cite{Ko90, Xu08}
from the discussion in Section 2. In particular,  $\varpi\equiv 0$ if $F$ is quasilinear. In this case,
(\ref{k-est}) becomes
\begin{equation}\label{k-est00} \k^c \ge
\min\{\k^0e^{Ac}; \k^1e^{A(c-1)}\}, \quad \forall \;c\in
[0,1].\end{equation} From the proof above, the strong maximum principle concludes that
if $"="$ holds for some $c_0\in (0,1)$ in (\ref{k-est00}), then $\kappa_s(x)\equiv constant$ for
all $x\in \Sigma^c$ and $\forall c\in (0,1)$. This implies that $\Sigma^{u(x)}$ is a
round sphere for every $x\in \bar{\Omega}$. The same conclusion is also true if condition (\ref{cdA-0})
is held. Note that $\varpi$ was used only in (\ref{new-000}) to get (\ref{h}). It's proved in \cite{BGMX} that $S_j\le 0$ under condition (\ref{cdA-0}). In that case, one may take $\varpi=0$ in (\ref{new-000}).

\bigskip

\noindent {\bf Acknowledgment:} We would like to thank Yong Huang for many helpful discussions.
Part of the work was done while
the second named author was visiting the McGill
University in 2009, she would like to thank the institution for the
warm hospitality.

\bibliographystyle{amsplain}

\end{document}